\documentclass[11pt,reqno]{amsart}

\usepackage[a4paper,inner=2.5cm,outer=2.5cm,top=2.5cm,bottom=2.5cm]{geometry}
\usepackage{amsmath,amssymb,amsthm,enumerate,mathtools,stmaryrd,xcolor}

\usepackage[pagebackref]{hyperref}
\hypersetup{colorlinks=true,linkcolor=blue,citecolor=teal,filecolor=magenta,urlcolor=cyan}

\usepackage[mathscr]{eucal}
\usepackage{bbm}

\usepackage{dsfont}

\usepackage[initials,nobysame]{amsrefs}

\usepackage{mathtools}
\usepackage[normalem]{ulem}

\usepackage{varioref}

\newcommand{\N}{\mathbb N}

\newcommand{\bw}{\mathbf{w}}

\newcommand{\Dim}{\rm Dim}

\theoremstyle{plain}
\newtheorem{thm}{Theorem}[section]
\newtheorem{lem}[thm]{Lemma}

\newtheorem{cor}[thm]{Corollary}

\theoremstyle{definition}

\theoremstyle{remark}

\numberwithin{equation}{section}

\usepackage[normalem]{ulem}

\date{} 
\begin{document}

\frenchspacing

\setlength\marginparsep{8mm}
\setlength\marginparwidth{20mm}

\title[Lyapunov exponents of matrix products]{Computation of Lyapunov exponents of matrix products
}

\author{Ai Hua Fan}
\address{(A. H. Fan) 
	LAMFA, UMR 7352 CNRS, University of Picardie, 33 rue Saint Leu, 80039 Amiens, France
    and
    Wuhan Institute for Math \& AI, Wuhan University, Wuhan 430072,  China}
\email{ai-hua.fan@u-picardie.fr}

\author{Evgeny Verbitskiy}
\address{(E. Verbitskiy)
	Mathematical Institute, Leiden University, Postbus 9512, 2300 RA Leiden, The Netherlands \& 
Korteweg-de Vries Institute for Mathematics, University of Amsterdam, Postbus 94248, 1090 GE
Amsterdam, The Netherlands
	}
\email{evgeny@math.leidenuniv.nl}

	\begin{abstract} 
		For $m$ given square matrices $A_0, A_1, \cdots, A_{m-1}$ ($m\ge 2$),
		one of which is assumed to be of rank $1$, and for a given sequence $(\omega_n)$ in $\{0,1, \cdots, m-1\}^\mathbb{N}$, the following limit, if it exists,  
$$L(\omega):=\lim_{n\to \infty} \frac 1n \log \|A_{\omega_0} A_{\omega_2}\cdots A_{\omega_{n-1}}\|$$ 
defines the Lyapunov  exponent of the sequence of matrices $(A_{\omega_n})_{n\ge 0}$. It is proved that the Lyapunov exponent $L(\omega)$ has a closed-form expression under certain conditions. One special case arises when $A_j$'s are non-negative and $\omega$ is generic with respect to some shift-invariant measure; a second special case occurs when $A_j$'s (for $1\le j<m$) are invertible and $\omega$ is a typical point  
 with respect to some shift-ergodic measure. 
Substitutive sequences and characteristic sequences of $\mathcal{B}$-free integers are  considered as examples. An application is presented for the computation of multifractal spectrum of weighted Birkhoff averages. 
	
		\end{abstract}

    \maketitle

	\section{Introduction}
 We study the Lyapunov exponent of products of  matrices in the standard setup of dynamical systems.
 	Let $(\Omega, T)$ be a topological dynamical system and let $A: \Omega \to M_{d} (\mathbb{C})$ be a continuous function from $\Omega$ into the space $M_{d} (\mathbb{C})$ of 
	$d\times d$ complex matrices. The (maximal)  Lyapunov exponent of $A$ at $\omega\in \Omega$ is defined by the limit 
	$$
	    L(\omega):=\lim_{n\to\infty} \frac{1}{n} \log \|A(\omega) A(T\omega) \cdots A(T^{n-1} \omega)\|,
	$$
	provided the limit exists.
	For any $T$-ergodic probability measure $\nu$, Furstenberg and Kesten \cite{FK1960} proved that	if $A$ takes values in the space ${\rm GL}(d,\mathbb{R})$  of invertible matrices and $\log^+\|A\| \in L^1(\nu)$, then the limit defining $L(\omega)$  exists 
	$\nu$-almost surely,  and is equal to the constant
	 $$
	     L_A:= \lim_n \frac{1}{n}  \mathbb{E}_\nu \log \|A^{(n)}(\omega)\|= \inf_n \frac{1}{n} \mathbb{E}_\nu  \log \|A^{(n)}(\omega)\|
	 $$   
	 where $$A^{(n)}(\omega) = A(\omega) A(T\omega) \cdots A(T^{n-1} \omega).$$	
	A relevant question is how to compute $L_A$. Another important question is for which `individual' points $\omega$ does the conclusion of the Furstenberg-Kesten  theorem hold. Does it holds for a given generic point for some invariant measure?
	Recall that a point  $\omega\in\Omega$ is said to be $\nu$-generic for some $T$-invariant measure $\nu$, if 
	$$
	      \frac{1}{n}\sum_{k=0}^{n-1} \delta_{T^k \omega} \to \nu
	$$
	in the weak-* topology. It is clear that 
	if  $\nu$ is ergodic, then
	$\nu$-almost every point $\omega$ is $\nu$-generic by the Birkhoff ergodic theorem. 
	
	Surprisingly, the conclusion of the Furstenberg-Kesten theorem is not necessarily true for a given
	$\nu$-generic point. 	 In fact,  even if $(\Omega, T)$ is uniquely ergodic (equivalently, if every point in $\Omega$ is $\nu$-generic), it is possible that the limit defining $L(\omega)$ fails to exist for some  $\omega$. Such counter-examples were first constructed by M. Herman \cite{Herman1981} and  P. Walters \cite{Walters1986}. 
	 
	For uniquely ergodic system, A. Furman \cite{Furman1997} identified sufficient conditions for
	uniform convergence of $n^{-1}\log\|   A( \omega)
	A( T\omega)\cdots A(T^{n-1}\omega)\|$ to the Lyapunov exponent $L_A$.  
    Recently, an affirmative answer regarding the existence of the Lyapunov exponent for every generic point was provided by Fan and Wu \cite{FW2022} in the case of non-negative matrices.
	
\medskip

	In general, computing Lyapunov exponents is a highly non-trivial task. As one of the problems, Kingman commented at the end of his paper \cite{Kingman1973} in 1973: ``\emph{Pride of place among the unsolved problems of subaddi­tive ergodic theory must go to the calculation of the Lyapunov exponent. In none of the applica­tions described here is there an obvious mechanism for obtaining an exact numerical value, and indeed this usually seems to be a problem of some depth}". 
	\medskip
	
	In this paper, we present a special class of products of matrices for which we can derive an explicit formula for the Lyapunov exponent. 
	This type of result is interesting for several reasons.
In particular, it can serve as a benchmark for numerical procedures to estimate the Lyapunov exponents.
\medskip

 Our sequence of matrices is selected from a finite set of matrices, one of which is assumed to be of rank one.  More precisely, let $\{A_0,A_1, \cdots, A_{m-1}\}$ be a set of
 $d\times d$
 matrices. We will assume  that $A_0$ is of rank $1$, namely
 $A_0= \mathbf{u}\mathbf{v}'\not=0$ where $\mathbf{u}$ and $\mathbf{v}$ are two  column vectors ($\mathbf{v}'$ denoting the transpose of $\mathbf{v}$).
Given $\omega =(\omega_k)_{k\ge 0} \in \{0,1, \cdots, m-1\}^\N$, we consider the corresponding sequence of matrices $(A_{\omega_k})_{k\ge 0}$.
 We will assume that the symbol $0$ appears infinitely often in $\omega$.  Our objective is to compute the Lyapunov exponent 
 \begin{equation}\label{eq:DefLyap}
	    L(\omega):=\lim_{n\to\infty} \frac{1}{n} \log \|A_{\omega_0} A_{\omega_1} \cdots A_{ \omega_{n-1}}\|.
	\end{equation}

Before stating our results, let us fix some notation and notions.  
We refer to $\{0,1, \cdots, m-1\}$ as the alphabet and let $\Omega=\{0,1,\ldots,m-1\}^\N$.
The left shift map $\sigma:\Omega\to\Omega$,
 $\omega'=\sigma(\omega)$, is given by $\omega'_n= \omega_{n+1}$ for all $n\in\N$. 
For a  word $\mathbf{w}=w_1\cdots w_\ell$ with $w_j\in\{0,\ldots,m-1\}$, we use $|\mathbf{w}|$ to denote the length $\ell$ of $\mathbf{w}$ and $A_{\mathbf{w}}$ to denote the matrix  product $A_{w_1}\cdots A_{w_\ell}$. The symbol $a^n$ will denote the word with the letter or word $a$,  concatenating  the letter or word $a$, n times.
We use $\mathcal{L}(\omega)$ to denote the language of $\omega$: namely, the set of all finite subwords contained in $\omega$.   If $\mathbf{w}$ is a proper subword of $\mathbf{w}'$, we write $\mathbf{w}  \prec \mathbf{w}'$.
For a given word $\mathbf{w}$, $[\mathbf{w}]$ denotes the cylinder  in $\Omega=\{0,\ldots,m-1\}^\N$, consisting of sequences having $\mathbf{w}$ as prefix. 
We use $1_{[\mathbf{w}]}$ to denote the indicator function of $[\mathbf{w}]$.

Consider $\omega\in\Omega$ such that the symbol $0$ appears infinitely often in $\omega$. Then $\omega$ can be uniquely written as  \begin{equation}\label{eq:omega}
 \omega = \omega_0\cdots \omega_{n-1} \cdots =
 \mathbf{w}_0 \, 0^{t_0}\, \mathbf{w}_1\, 0^{t_1}\mathbf{w}_2\, 0^{t_2}\cdots \mathbf{w}_{r}\, 0^{t_{r}}\cdots 
\end{equation}
where the words $\{\mathbf{w}_j\}$ do not contain $0$ (with $\mathbf{w_0}$ possibly an empty word) and 
$\{t_j\}$ are positive integers. Let $\mathcal R=\mathcal R(\omega)=\{\mathbf w_1,\mathbf w_2,\ldots\}$. The words in $\mathcal{R}$ will be called return words (to $0$).
\medskip

When $A_0, A_1, \cdots, A_{m-1}$ are non-negative matrices
and $\omega$ is a generic point for some ergodic measure,  the following theorem is proved.

 \begin{thm} \label{thm:positive}

 Suppose that 

 \begin{itemize}
 \item[(1)]  $\{A_0,A_1, \cdots, A_{m-1}\}$ is a collection of non-negative matrices; 
$A_0= \mathbf{u}\mathbf{v}'$ for some vectors $\mathbf{u},\mathbf{v}$ such that $\mathbf{v}'\mathbf{u}\not=0$,
\item[(2)]  $\omega=\{\omega_n\}_{n\ge 0}$ is  a point which is generic  for some ergodic measure $\nu$ with $\nu([0])>0$,
\end{itemize}
then  the limit \eqref{eq:DefLyap} defining $L(\omega)$ exists. 
Moreover, $L(\omega)=-\infty$ if 
 $A_{{\mathbf w}_0}\mathbf{u}=0$ or
$\mathbf{v}'A_{\mathbf{w}} \mathbf{u}=0$ for some
$\mathbf{w}\in\mathcal R$;  
otherwise, $L(\omega)$ is given by
\begin{equation}\label{eq:LiapRankMain}
  L (\omega)= 
  \nu([0]) \log \mathbf{v}'\mathbf{u} +\sum_{\mathbf{w}\in \mathcal{R}}  \nu([0\mathbf{w}0])\log \frac{\mathbf{v}' A_{\mathbf{w}} \mathbf{u}}
   {\mathbf{v}'\mathbf{u} } \in [-\infty, +\infty).
\end{equation}
\end{thm}

If we assume that $\{A_1, \cdots, A_{m-1}\}$ are all invertible,
we have the following formula for the Lyapunov exponent of
products of random matrices. 

\begin{thm}\label{thm:random}
    Let $\{A_0,A_1, \cdots, A_{m-1}\}$ be a set of
 complex  matrices. Suppose that  
$A_0= \mathbf{u}\mathbf{v}'$ is of rank one such  that $\mathbf{v}'\mathbf{u}\not=0$ and that $A_1, \cdots, A_{m-1}$
are invertible. Let $\mu$ be a shift invariant and ergodic measure on $\{0, 1, \cdots, m-1\}^\N$ such that 
$\mu([0])>0$. Then for $\mu$-almost every $\omega$,  the Lyapunov exponent $L(\omega)$ defined by (\ref{eq:DefLyap})
exists and is equal to
\begin{equation}\label{eq:invertible}
  L = 
  \mu([0])  \log |\mathbf{v}'\mathbf{u}| +\sum_{\mathbf{w} \in \mathcal{R}}   \mu([0\mathbf{w}0]) \log \frac{|\mathbf{v}' A_{\mathbf{w}} \mathbf{u}|}
   {|\mathbf{v}'\mathbf{u} |} \in [-\infty, +\infty),
\end{equation}
where $\mathcal{R}$  is the set of  all words $\mathbf{w}$  not containing $0$ such that $\mu([0\mathbf{w}0])>0$.

\end{thm}

Both Theorem \ref{thm:positive} and Theorem \ref{thm:random} are corollaries of a more general result  (Theorem \ref{thm:rank1}). 
\medskip

The following result is a corollary of Theorem \ref{thm:random}.

\begin{cor} Under the assumption of Theorem \ref{thm:random}, if, furthermore, $\mu$
is the  Bernoulli measure corresponding to the probability vector $(p_0, p_1, \cdots, p_{m-1})$ with  $p_0 >0$.
Then almost surely we have
\begin{equation}\label{eq:LiapRank-Bernoulli}
  L (\omega)= 
  p_0 \log |\mathbf{v}'\mathbf{u}| +\sum_{\mathbf{w}\in \mathcal{R}}  p_0^2 p_{\mathbf{w}} \log \frac{|\mathbf{v}' A_{\mathbf{w}} \mathbf{u}|}
   {|\mathbf{v}'\mathbf{u} |}
\end{equation}
where $p_{\mathbf{w}} = p_{w_1}p_{w_2}\cdots p_{w_\ell}$ for $\mathbf{w} = w_1w_2\cdots w_\ell$.
\end{cor}

It is also easy to obtain   similar expressions for Lyapunov exponents in case of  Markov measures $\mu$ with $\mu([0])>0$.

\medskip

Expressions for Lyapunov exponents, similar to (\ref{eq:invertible})
and (\ref{eq:LiapRank-Bernoulli}), 
have been obtained in a number of particular cases of random matrix products, where one or all
matrices are presumed to be of rank one. For example, Pincus \cite{Pincus} (1985)
proved that the Lyapunov exponent
for random products of two $2\times 2$ matrices
$$
A_0=\begin{bmatrix} 1 & 0\\ 0 & 0\end{bmatrix},
\quad A_{1}=\text{arbitrary}, 
$$
chosen independently with $\mathbb P(A_0)=p$, $\mathbb P(A_1)=1-p$, 
is given by
$$
L=\sum_{n=1}^\infty p^2(1-p)^{n} \log|  {\mathbf u}' A_1^n  {\mathbf u}|,\quad{\mathbf u}=\begin{pmatrix}1\\ 0\end{pmatrix}.
$$
This formula coincides with \eqref{eq:LiapRank-Bernoulli}. 
R. Lima and M. Rahibe \cite{LIMA} (1994)  generalized this result
to the case when matrices $A_0,A_1$ are chosen according to some Markov measure. The formula \eqref{eq:invertible}
covers the case of Markov measures.
E. Key \cite{KEY} (1987) considered independent random products of two  matrices $A_0$, $A_1\in M_{d\times d}(\mathbb C)$, satisfying the following algebraic relations: 
$$A_0^2=zA_0\text{ and }A_0A_1^kA_0=x_kA_0A_1A_0+y_kA_0
$$ 
for some (complex) numbers $\{z,x_k,y_k| k\ge 1\}$ and obtained
$$
L=\lim _{n \rightarrow \infty} n^{-1} \log \left\|A_{\omega_n} \ldots A_{\omega_1}\right\|=
\max_{s \in \mathbb S^1}(1-p)^2 \sum_{k=0}^{\infty} \log \left|x_k s z^{-1}+y_k\right| p^k.
$$
Notice that if $A_0$ has rank 1, many matrices $A_1$ satisfy the above algebraic conditions.
More recently, Altschuler and Parrilo \cite{Alts} (2020) considered independent products of any number of rank one matrices $A_j =  {\mathbf u}_j {\mathbf v}_j'$ ($0\le j<m)$,
and found that the corresponding Lyapunov exponent is given by
$$
L =\sum_{i,j=0}^{m-1} p_ip_j \log\bigm|   {\mathbf u}_i' {\mathbf v}_j \bigm|, 
$$
where $p_j$ is the probability to choose the matrix $A_j$, $j=0,\ldots, m-1$.
We point out that Pollicott \cite{Pollicott2010} (2010) developed  an 
accurate numerical method for computing the Lyapunov exponent for products of independent identically distributed random matrices. 
 \medskip

 Let us look at the application of Theorem \ref{thm:random} 
to products of random matrices arising in quantum mechanics
\cite{Tristan}. Consider the set of two matrices 
 
 $$
     A_0 =\begin{pmatrix}
         1 & 0 \\ 0 & 0
     \end{pmatrix}, \qquad
   A_1 = R_\theta = \begin{pmatrix}
         \cos \theta  & -\sin \theta \\ \sin \theta & \ \ \ \cos \theta
     \end{pmatrix}.
     $$

     We consider the invariant Markov measure $\mu$ defined by the transition matrix
     $$
         P =\begin{pmatrix}
         p_0 & 1-p_0 \\ 1-p_1 & p_1
     \end{pmatrix} \quad (0<p_0, p_1 <1).
     $$
     The unique invariant probability vector $\pi=(\pi_0, \pi_1)$,
     which satisfies $\pi P=\pi$, is such that
     $\pi_0 =\frac{1-p_1}{2-(p_0+ p_1)}$ and $\pi_1 =\frac{1-p_0}{2-(p_0+ p_1)}$. 
     When $1-p_0=p_1=p$ for some $0<p<1$, the resulting measure is Bernoulli. By Theorem \ref{thm:random}, for $\mu$-almost all $\omega$, the Lyapunov exponent $L(\omega)$ is equal to 
     $$
          \lambda(\theta) = \sum_{n=1}^\infty \pi_0 (1-p_0) p_1^n \log \frac{1}{|\cos n \theta|}.
     $$
     Notice that $\lambda(\theta) \in [-\infty, 0)$ for all $\theta$, meaning that the matrix product is contractive, even ``strongly" contractive (i.e. $\lambda(\theta) =-\infty$). Moreover, for the Lebesgue almost all
     $\theta$, $\lambda(\theta)$ is finite because $\log \frac{1}{|\cos \theta|}$ is integrable on the interval $[0, 2\pi]$. On the other hand, it is easy to see that there exists an uncountable dense set of $\theta$ such that $\lambda(\theta)=-\infty$. For example, we can choose a  sequence of integers $(n_k)$ that is highly lacunary in the sense that 
     $$
            n_{k+1} \gg n_1 + \cdots + n_k
     $$
     so that $p_1^{2^{n_k}} \log \frac{1}{|\cos 2^{n_k} \theta|}\ge 1$
     for $\theta$ defined by $$\frac{\theta}{\pi} = \sum_{k=1}^\infty \frac{1}{2^{n_k}}.$$
     Consequently the function $\lambda(\cdot)$ is nowhere continuous.
     \medskip 

 Theorem   \ref{thm:positive} -- computing the Lyapunov exponent for a ``deterministic" matrix product, seems entirely new. Here by ``deterministic"
we mean  a fixed sequence of matrices, while ``random"  means an almost sure sampling of matrices. 
\medskip 

 Theorem   \ref{thm:positive} and Theorem \ref{thm:random} can be applied in many interesting cases. For example, as we will see  in Section  \ref{sect:SS} and Section
 \ref{sect:BF}, respectively,  when $\omega$ is a primitive substitutive sequence or when   $\omega$ is the characteristic function of the set of $\mathcal{B}$-free integers.  Section \ref{sect:MA}  is devoted to an application to the computation of 
multifractal spectrum of weighted Birkhoff averages.  
In Section \ref{Sect:GF}, we will first  prove a general result (Theorem \ref{thm:rank1})
which covers Theorem   \ref{thm:positive} and Theorem \ref{thm:random}. The proofs of Theorem   \ref{thm:positive} and Theorem \ref{thm:random} will be given in Section
\ref{Sect: Proof1} and Section  \ref{Sect: Proof2}. Section
\ref{Sect: ExactF} provides a way to compute the measures
$\nu([0\mathbf{w}0])$, called exact frequency of the return word $\mathbf{w}$, which appear in our formulas. 
\medskip

{\em Acknowledgement.}  We thank Fabien Durand for pointing out the method of derivative substitution for computing the exact frequencies. The first author is partially supported by NSFC (no. 12231013).

\section{A general formula} \label{Sect:GF}

In this section, we prove a general formula for Lyapunov exponent. Our proof is based on the analysis of return words.

\subsection{Decomposition of \texorpdfstring{$\omega$}{omega} into return words and return word statistics.}

Let $\omega \in \{0,1, \cdots, m-1\}^\N$ be a sequence in which $0$ appears infinitely often. Consider the decomposition  \eqref{eq:omega} which is repeated below:
\begin{equation}\label{eq:decomp2}
 \omega = \omega_0\cdots \omega_n \cdots =
 \mathbf{w}_0 \, 0^{t_0}\, \mathbf{w}_1\, 0^{t_1}\mathbf{w}_2\, 0^{t_2}\cdots \mathbf{w}_{r}\, 0^{t_{r}}\cdots 
\end{equation}
where the words $\{\mathbf{w}_j\}$ do not contain $0$ (with $\mathbf{w_0}$ possibly an empty word) and 
$\{t_j\}$ are positive integers. We define the set of return word $\mathcal R$ as the set of all $\bw$-words which appear in the decomposition
above, but not including the first word $\bw_0$. Equivalently, if we define  the language of 
$\omega$, denoted by $\mathcal{L}(\omega)$,  as the collection of all finite subwords appearing in $\omega$, then
 $$
    \mathcal{R} =\{\mathbf{w}\in \mathcal{L}(\omega): 0\mathbf{w} 0\in  \mathcal{L}(\omega)\quad \text{and} \quad\mathbf{w} \ {\rm does \ not \ contain\ } 0 \}.
$$
Let $\mathcal{T}$ be the set of subwords of $\omega$ of the form $0^t$, which appear between two consecutive return words, i.e.,
 $$
    \mathcal{T} =\{0^t\in \mathcal{L}(\omega): \mathbf{w}0^t\mathbf{\tilde w} \in  \mathcal{L}(\omega) \ {\rm for \ some} \ \mathbf{w}, \mathbf{\tilde w} \in \mathcal{R} \}.
$$
Note that if $\omega$ is minimal, both sets $\mathcal{R}$ and  $\mathcal{T}$ are finite.   In this case, the sum in \eqref{eq:LiapRankMain} is taken over a finite set.

\medskip
We define the (asymptotic) frequency $\rho_{\mathbf a}$ of a finite word $\mathbf a$  in $\mathcal{L}(\omega)$ as follows:
$$\aligned
  \rho_{\mathbf a} &:=
   \lim_{n\to \infty}\frac{1}{n} N_n(\mathbf a|\omega)
   =\lim_{n\to\infty}\frac{1}{n}\sum_{j=0}^{n-1} 1_{[\mathbf{a}]} (\sigma^j \omega),
   \endaligned
$$
provided the limit exists, 
where
$$
N_n(\mathbf a|\omega) :=\#\bigl\{0\le j<n-1:  \mathbf{a} \ \text{is a prefix of }\sigma^j\omega\bigr\}.
$$

   For any return word $\bw$  in $\mathcal{R}$,  we will need a different frequency: namely, the frequency of appearances of the return word $\bw$ between two consecutive symbols $0$ in the decomposition \eqref{eq:decomp2},
   which will be referred to  as  the exact frequency of $\bw\in\mathcal R$ and will be denoted by $F_\mathbf{w}$. It is clear that  $F_\mathbf{w}$ is  the
   asymptotic frequency of a word $0\bw 0$:
\begin{equation}\label{eq:F=rho}
       F_{\mathbf{w}}  = \rho_{0\mathbf{w}0}. 
\end{equation}
The exact frequency $F_\mathbf{w}$  does not account for occurrences of $\mathbf{w}$  as a proper subword of some longer return words $\mathbf{w}'\in \mathcal{R}$. 

If $\omega$ is a generic point for some translation invariant measure $\nu$, we have
\begin{equation}\label{eq:Fw-nu}
\rho_{0}=\nu([0]), \qquad 
F_{\mathbf w}=\nu([0\bw 0]). 
\end{equation}
 
\subsection{Main result}

The following theorem describes conditions under which we are able to derive  explicit expressions
 for  Lyapunov exponents.
 Its proof is based on the decomposition (\ref{eq:omega}) of $\omega$.  
For a large integer, say $n\ge |\mathbf{w}_0| + {t_0}\ge 1$, we can write the $n$-prefix of $\omega$ as follows
\begin{equation}\label{eq:omega-n}
 \omega|_0^{n-1}=\omega_0\cdots \omega_{n-1} =\mathbf{w}_00^{t_0}\mathbf{w}_10^{t_1}\mathbf{w}_20^{t_2}\cdots \mathbf{w}_{r}0^{t_{r}}\mathbf{w}^*_{r+1} 0^{t^*_{r+1}}
\end{equation}
where $r=r(n)$ depends on $n$,  $0^{t_j} \in \mathcal{T}$  for $0\le j\le r$ (therefore $t_j\ge 1$) and $\mathbf{w}_j \in \mathcal{R}$ for $1\le j\le r$,   but $t^*_{r+1}\ge 0$
and $\mathbf{w}_{r+1}^*$ is a prefix of a return word (it is a return word when $t_{j+1}^*\ge 1$).
  
  \begin{thm} \label{thm:rank1}

 Let $\{A_0,A_1, \cdots, A_{m-1}\}$ be a set of
complex  matrices. Suppose that  
$A_0= \mathbf{u}\mathbf{v}'$ is of rank one such that the inner product $\mathbf{v}'\mathbf{u}\not=0$.
Let  $\omega=\{\omega_n\}_{n\ge 0}$ is a point of $\Omega=\{0,1, \cdots, m-1\}^\N$ 
such that $0$ appears infinitely many times in  $\omega$
with $(\mathcal R,\mathcal T$)-decomposition  (\ref{eq:decomp2}) such that the frequency $\rho_0$ exists 
and the exact frequencies $F_{\mathbf{w}}$ exist for 
all $\mathbf{w}\in \mathcal{R}$.

Suppose furthermore that
\begin{enumerate}
\item[{\rm (a)}] 
 $A_{\mathbf{w_0}}\mathbf{u}\not=0$; 
$\mathbf{v}'A_{\mathbf{w}_j}
\mathbf{u}\not=0$ for all $j\ge 1$;
\item[{\rm (b)}] $\gamma_n:=\max_{\mathbf{w}' \prec\mathbf{w}_{r+1}} |\log  \|  \mathbf{v}'A_{\mathbf{w}'}\||=o(n)$ as $n\to \infty$;
\item [{\rm (c)}] 
we have
$$
    \lim_{M\to \infty} \varlimsup_{n\to \infty} \frac{1}{n}\sum_{\mathbf{w}_j \prec \omega|_0^{n-1}, |\mathbf{w}_j|>M} |\log |\mathbf{v}'A_{\mathbf{w}_j} \mathbf{u}|| =0.
$$
\end{enumerate}
Then  the following limit exists
$$
    L= \lim_{n\to \infty} \frac{1}{n} \log \|A_{\omega_0}A_{\omega_2}\cdots A_{\omega_{n-1}}\|  \in [-\infty, +\infty)
$$
and is given by 
\begin{equation}\label{eq:LiapRank1}
  L = 
  \rho_0 \log |\mathbf{v}'\mathbf{u}| +\sum_{\mathbf{w}\in \mathcal{R}}  F_{\mathbf{w}} \log \frac{|\mathbf{v}' A_{\mathbf{w}} \mathbf{u}|}
   {|\mathbf{v}'\mathbf{u} |}.
\end{equation}
\end{thm}

 {\bf Remark 1.} We have trivially $L=-\infty$ if $A_{\bw_0}\mathbf{u} =0$ or $\mathbf{v}'A_{\bw_j}\mathbf{u} =0$ for some $j$ (i.e. the condition (a) is violated), or $\mathbf{v}'\mathbf{u}=0$. Because, as we shall see in the proof of Theorem \ref{thm:rank1},  
 for large $n$ we have
$$
    A_{\omega_0} A_{\omega_2}\cdots  A_{\omega_n-1}=0.
$$

{\bf Remark 2.}  
Assume that the set $\mathcal{R}$ is finite. Then both  the conditions (b) and  (c) are trivially satisfied,   the expression (\ref{eq:LiapRank1}) involves only a finite number of matrix products $A_{\mathbf{w}}$, which can be directly evaluated, and a finite number of exact frequencies $F_{\mathbf{w}}$ is to be computed. Therefore, when $\mathcal{R}$ is finite,  the formula \eqref{eq:LiapRank1} is effective. 
Thus, Theorem \ref{thm:rank1} applies to all minimal sequences, 
including primitive substitutive sequences $\{\omega\}$, like Thue-Morse sequence and Fibonacci sequence.

{\bf Remark 3.}  
 As we shall see in  
  the proof of Theorem \ref{thm:rank1},
  the positive part of the  series \eqref{eq:LiapRank1} converges. Recall that for a real series $\sum a_n$ with $a_n \in \mathbb{R}$, its positive part  is  $\sum a_n^+$  where $x^+$ 
denotes   the positive  part of a real number $x$.

{\bf Remark 4.} 
The closed formula (\ref{eq:LiapRank1})
 can be restated as follows
 \begin{equation}\label{eq:LiapRank2}
   L = 
  \left(1-  \sum_{\mathbf{w} \in \mathcal{R}}  (|\mathbf{w}|+1)F_{\mathbf{w}}\right)  \log \mathbf{v}'\mathbf{u} +\sum_{\mathbf{w}\in \mathcal{R}}  F_{\mathbf{w}} \log \mathbf{v}' A_{\mathbf{w}} \mathbf{u}.
\end{equation}

\subsection{A counterexample} \label{sect:counterexample}

 In general, the formula \eqref{eq:LiapRank1} does not hold when the condition (c) in Theorem \ref{thm:rank1} is violated, even for positive matrices.

Consider the  following sequence
\begin{equation}\label{eq:omegaa}
  \omega=0101101110\cdots 01^n01^{n+1}\cdots.
\end{equation}
It is clear that  $\mathcal{R}=\{1, 11, 111, \cdots\}$, $F_\mathbf{w}=0$ for all $\mathbf{w}\in \mathcal{R}$ and $\rho_0=0$.
The following example shows that the formula \eqref{eq:LiapRank1} in Theorem \ref{thm:rank1} is not always true
when the condition (c) in Theorem \ref{thm:rank1} is violated.
\medskip

{\bf Example.}   Consider the above sequence $\omega$ defined by \eqref{eq:omegaa}.
Consider two $d\times d$ matrices: $A_0 = \mathbf{1} \mathbf{1}'$, and $A_1$ being a nonnegative and invertible matrix. Then
\begin{equation}\label{eq:example1}
    L(\omega) = \rho(A_1),
\end{equation}
where $\rho(A_1)$ is the spectral radius of $A_1$.  But the formula  \eqref{eq:LiapRank1} would give $L(\omega)=0$.
\medskip

Indeed, if we consider the product of  first $\frac{n(n+1)}{2} +n+1$ matrices of $(A_{\omega_j})$, we get
\begin{eqnarray*}
\|A_{ 0} A_{ 1}A_{ 0}A_{1}^2A_{ 0}A_{ 1}^3A_{\ 0} \cdots A_{ 0}A_{ 1}^nA_{ 0}\| 
& = & \mathbf{1} (\mathbf{1}' A_1 \mathbf{1})  (\mathbf{1}' A_1^2 \mathbf{1}) \cdots  (\mathbf{1}' A_1^n \mathbf{1}) \mathbf{1}'\\
&=&\|A_1\| \|A_1^2\|\cdots \|A_1^n\| \|\mathbf{1} \mathbf{1}'\|
\end{eqnarray*}
where $\|A\|$ is chosen to be  the norm $\sum_{i. j} |a_{i,j}|$ for a matrix $A=(a_{i,j})$.
Since $ \|\mathbf{1} \mathbf{1}'\|=d^2$ is a positive constant and $\log \|A_1^n\| \sim n \rho(A_1)$ as $n \to \infty$, we get 
\begin{equation}\label{eq:A1}
 \lim_{n\to \infty} \frac{ \log (\|A_1\| \|A_1^2\|\cdots \|A_1^n\| \|\mathbf{1} \mathbf{1}'\|)}{\frac{n(n+1)}{2} +n+1}
   = \rho(A_1).
\end{equation}
Consider  $1\le m\le n+1$ and the product
$$
\|A_{ 0} A_{ 1}A_{ 0}A_{1}^2A_{ 0}A_{ 1}^3A_{\ 0} \cdots A_{ 0}A_{ 1}^nA_{ 0} A_1^m\| 
=\|A_1\| \|A_1^2\|\cdots \|A_1^n\| \|\mathbf{1} \mathbf{1}' A^m\|.
$$
With the help of \eqref{eq:A1}, in order to conclude for \eqref{eq:example1}, we only need to check that 
$$
    \lim_{n\to \infty} \frac{1}{n^2} \max_{1\le m \le n+1}|\log \|\mathbf{1} \mathbf{1}' A_1^m\||=0.
$$
The checking is as follows. 
On the one hand, evidently we have
$$
    \|\mathbf{1} \mathbf{1}' A_1^m\| \le d^2 \|A_1\|^m,
$$
and on the other hand, as $A_1$ is invertible, we have $ \|\mathbf{1} \mathbf{1}' \|\le  \|\mathbf{1} \mathbf{1}' A_1^m\| \|A_1^{-m}\|$ and then
$$
             \|\mathbf{1} \mathbf{1}' A_1^m\| \ge \frac{d^2} {\|A_1^{-1}\|^m}.
$$
It follows that $\log   \|\mathbf{1} \mathbf{1}' A_1^m\| =O(m) =o(n^2)$.
This last estimate actually shows that the condition (b) is satisfied. 

Assume furthermore that $A_1$ is non negative. Then
$\log |\mathbf{1}' A_1^m \mathbf{1}|\sim \rho(A_1)$ as $m\to \infty$. 
It can be deduced that the limit involved in the condition (c) exists and is equal to $\log \rho(A_0)$.
Thus the condition (c) of Theorem \ref{thm:rank1} is not satisfied if  $\rho(A_1) \not =1$.

 \subsection{Convergence of the  series $\sum_{\mathbf{w}\in \mathcal{R}}F_\mathbf{w} \log \left|\frac{\mathbf{v}'A_{\mathbf{w}} \mathbf{u} }{\mathbf{v}' \mathbf{u}}\right|$
 } 
 
 Before proving Theorem \ref{thm:rank1}, let us first discuss the convergence of the series appearing in  Theorem \ref{thm:rank1}.
 As we shall see,  its positive part is convergent and thus its sum is finite or $-\infty$.

The following notation will be used several times.  For any integer $M\ge 1$, denote
$$
  \mathcal{R}_M =\{ \mathbf{w}\in\mathcal{R}: |\mathbf{w}|\le M\}.
$$
It is the set of  all return words of length not exceeding  $M$.

\begin{lem} \label{lem:lem1}$\sum_{\mathbf{w}\in \mathcal{R}} |\mathbf{w}| F_\mathbf{w} \le 1-\rho_0$.
\end{lem}

\begin{proof} Let us first observe that
\begin{equation}\label{eq:tw}
    |\mathbf{w}_0| + |\mathbf{w}_1| +\cdots + |\mathbf{w}_r| +|\mathbf{w}_{r+1}^*|    + t_0+t_1 +\cdots +t_r +t_{r+1}^* =n.
\end{equation}
Also notice that $N_n(0)= t_0+t_1 +\cdots +t_r +t_{r+1}^*$ is the number of $0$'s contained in $\omega|_0^{n-1}$ and 
$\mathbf{w}_j$'s (possibly $\mathbf{w}_{r+1}^*$) in \eqref{eq:tw} are all return words contained in $\omega|_0^{n-1}$.
We keep only those return words with length not exceeding $M$ to get 
     $$
        \frac{1}{n} \sum_{\mathbf{w}\in \mathcal{R}_M} |\mathbf{w}| N_n(\mathbf{w}) \le 1- \frac{N_n(0)}{n}.
     $$
     We get an equality if $\mathbf{w}_{r+1}^*$ is a return word. But if $\mathbf{w}_{r+1}^*$ is not a return word and is very long, we get a big gap in the last inequality. See Example in Section \ref{sect:counterexample}.
     Since it is assumed that the exact frequencies $F_\mathbf{w}$ exist, letting $n\to \infty$ we get  $\sum_{\mathbf{w}\in \mathcal{R}_M} |\mathbf{w}| F_\mathbf{w} \le 1- \rho_0$. 
     The claim is proved since $M$ is chosen arbitrarily.
\end{proof}

\begin{lem} \label{lem:lem2}
The series $\sum_{\mathbf{w}\in \mathcal{R}}F_\mathbf{w} \log \left|\frac{\mathbf{v}'A_{\mathbf{w}} \mathbf{u} }{\mathbf{v}' \mathbf{u}}\right|$ has
its sum  in the interval $[-\infty, +\infty)$.
We can state the result as follows: the following limit exists
$$
    \lim_{M\to \infty} \sum_{\mathbf{w}\in \mathcal{R}_M}F_\mathbf{w} \log \left|\frac{\mathbf{v}'A_{\mathbf{w}} \mathbf{u} }{\mathbf{v}' \mathbf{u}}\right| \in [-\infty, +\infty).
$$
\end{lem}

\begin{proof}  Since $|\mathbf{v}'A_{\mathbf{w}} \mathbf{u} | \le \|\mathbf{v}\|  \|\mathbf{u}\| \Delta^{|\mathbf{w}|}$ where 
$\Delta =\max_{1\le j\le m-1}\|A_j\|$,  there exists a constant $D>0$ such that
    $$
        \log \left|\frac{\mathbf{v}'A_{\mathbf{w}} \mathbf{u} }{\mathbf{v}' \mathbf{u}}\right|< D |\mathbf{w}|.
    $$
    Thus the following series of negative terms verifies 
    $$
         \sum_{\mathbf{w}\in \mathcal{R}}F_\mathbf{w}\left( \log \left|\frac{\mathbf{v}'A_{\mathbf{w}} \mathbf{u} }{\mathbf{v}' \mathbf{u}}\right| - D |\mathbf{w}| \right)
         \in [-\infty, 0).
    $$
    We finish the proof by using the fact that  the series $\sum_{w \in\mathcal{R}} |\mathbf{w}| F_{\mathbf{w}}$ converges (Lemma \ref{lem:lem1}). 
\end{proof}

\subsection{Proof of  Theorem \ref{thm:rank1}}

The definition of Lyapunov exponent does not depend on the choice of matrix norm. We will  choose the 
most convenient  norm for us, namely,  the Frobenius norm or the Hilbert-Schmidt norm  defined by 
$$
   \|A\|_F = \sqrt{{\rm tr}(A^*A)}
$$
where $A=(a_{ij})$ and $A^*=(\overline{a_{ji}})$. 
Note that for a rank one matrix $A=\mathbf{u} \mathbf{v}'$ one has
\begin{equation}\label{eq:FNab}
  \|\mathbf{u} \mathbf{v}'\|_F =\|\mathbf{u}\|_2 \|\mathbf{v}\|_2
\end{equation}
where $\|\cdot\|_2$ is the Euclidean norm.
Furthermore, for any $t\ge 1$  we have
\begin{equation}\label{eq:A_0^t}
    A_{0^t} = (\mathbf{u}\mathbf{v}')^t =  \mathbf{u} (\mathbf{v}' \mathbf{u})^{t-1} \mathbf{v}'= \lambda^{t-1} \mathbf{u} \cdot \mathbf{v}' =\lambda^{t-1} A_0
 \end{equation}
 where 
 $\lambda=\mathbf{v}'\mathbf{u}=\sum_{j=1}^d u_jv_j$ is a scalar. 
 
 Therefore, using the notation introduced in (\ref{eq:omega-n}) and the associativity of matrix multiplication, one has
 $$\aligned
 A_{\omega_0} A_{\omega_2}&\cdots A_{\omega_{n-1}}=A_{\mathbf{w}_0} A_0^{t_0}A_{\mathbf{w}_1} A_0^{t_1}\cdots
 A_{\mathbf{w}_{r}}A_0^{t_{r}}A_{\mathbf{w}^*_{r+1}0^{t^*_{r+1}}}\\
 &=A_{\mathbf{w}_0} (\lambda^{t_0-1} \mathbf{u}  \mathbf{v}') A_{\mathbf{w}_1} (\lambda^{t_1-1} \mathbf{u}  \mathbf{v}')\cdots A_{\mathbf{w}_{r}} (\lambda^{t_r-1} \mathbf{u}  \mathbf{v}')A_{\mathbf{w}^*_{r+1}0^{t^*_{r+1}}}\\
 &=\lambda^{\sum_{j=0}^r (t_j-1)} A_{\mathbf{w}_0}  (\mathbf{u}  \mathbf{v}') A_{\mathbf{w}_1}  (\mathbf{u}  \mathbf{v}')\cdots A_{\mathbf{w}_{r}}  (\mathbf{u}  \mathbf{v}')A_{\mathbf{w}^*_{r+1}0^{t^*_{r+1}}}\\
 &=\lambda^{\sum_{j=0}^r (t_j-1)} (A_{\mathbf{w}_0}\mathbf{u}) ( \mathbf{v}' A_{\mathbf{w}_1} \mathbf{u})
 \cdots  ( \mathbf{v}' A_{\mathbf{w}_r} \mathbf{u})  ( \mathbf{v}' A_{\mathbf{w}^*_{r+1}0^{t^*_{r+1}}}). \\
 \endaligned
 $$
 Taking into account that all matrix products $\mathbf{v}' A_{\mathbf{w}_j} \mathbf{u}$, $j=1,\ldots,r$, are scalars, we conclude that
 \begin{eqnarray*}
 \|A_{\omega_0}\cdots A_{\omega_{n-1}}\|_F
 &=& |\lambda|^{\sum_{j=0}^r (t_j-1)}\prod_{j=1}^r  | \mathbf{v}' A_{\mathbf{w}_j} \mathbf{u}| \times \| A_{\mathbf{w}_0}\mathbf{u}\mathbf{v}' A_{\mathbf{w}^*_{r+1}0^{t^*_{r+1}}}\|_F
 =: M_n\times R_n
 \end{eqnarray*}
 where 
 $$
 M_n :=
 |\lambda|^{\sum_{j=0}^r (t_j-1)}\prod_{j=1}^r  | \mathbf{v}' A_{\mathbf{w}_j} \mathbf{u}|, 
 \qquad R_n := \| A_{\mathbf{w}_0}A_0 A_{\mathbf{w}^*_{r+1}0^{t^*_{r+1}}}\|_F.
 $$
We will refer to $M_n$ as the main term, and to $R_n$ as the remainder term. We will show that the remainder has no contribution to the Lyapunov exponent.

 Let us start by analysing the remainder term $R_n$, which contains at least one matrix $A_0$, and hence is
 of rank at most one. It is actually of rank one i.e. $R_n>0$, as we see now. Two situation may occur. First, if $t_{r+1}^*=0$, one has 
$$
A_{\bw_0}A_0A_{\bw_{r+1}^* 0^{t_{r+1}^*}}  =
(A_{\bw_0} \mathbf{u})  ( \mathbf{v}' A_{\bw_{r+1}^*}),
$$ 
and thus by \eqref{eq:FNab} and the condition (a), 
$$
  R_n = \|A_{\bw_0}\mathbf{u}\|_2 \|\mathbf{v}'A_{\bw_{r+1}^*}\|_2 >0. 
$$
Secondly,  if $t_{r+1}^*>0$, then $\bw_{r+1}^*=\bw_{r+1}$ and 
$$\aligned
A_{\bw_0}A_0&A_{\bw_{r+1} 0^{t_{r+1}^*}}  =(A_{\bw_0}\mathbf{u})\cdot (\mathbf{v}'A_{\bw_{r+1}}) \cdot \lambda^{t_{r+1}^*-1} \mathbf{u}\mathbf{v}'\\
&=\lambda^{t_{r+1}^*-1} (A_{\bw_0}\mathbf{u})\cdot (\mathbf{v}'A_{\bw_{r+1}}  \mathbf{u}) \cdot\mathbf{v}
=\lambda^{t_{r+1}^*-1}(\mathbf{v}'A_{\bw_{r+1}}  \mathbf{u}) (A_{\bw_0}\mathbf{u})\mathbf{v}',
\endaligned
$$
since $\mathbf{v}'A_{\bw_{r+1}}  \mathbf{u}$ is a scalar. Thus, applying \eqref{eq:FNab} again, one has
$$
R_n =|\lambda|^{t_{r+1}^*-1}\cdot|\mathbf{v}'A_{\bw_{r+1}}  \mathbf{u}|\cdot 
\|A_{\bw_0}\mathbf{u}\|_2 \|\mathbf{v}\|_2 >0.
$$
Therefore, depending on whether the last symbol in $(\omega_0,\ldots, \omega_{n-1})$ is the symbol 0 or not, i.e.,
whether, $t_{r+1}^*>0$ or $t_{r+1}^*=0$, we meet two situations. First, if $t_{r+1}^*=0$,  
we  obtain that the equality
$$
  \|A_{\omega_0}\cdots A_{\omega_{n-1}}\|_F=
 |\lambda|^{\sum_{j=0}^r (t_j-1)}\Bigl( \prod_{j=1}^r  | \mathbf{v}' A_{\mathbf{w}_j} \mathbf{u}| \Bigr)\cdot 
\|A_{\bw_0}\mathbf{u}\|_2\cdot \|\mathbf{v}'A_{\bw_{r+1}^*}\|_2
$$
where $\lambda=\mathbf{v}'  \mathbf{u}$, 
that we write in the following way
\begin{equation}\label{expression1}
\|A_{\omega_0}\cdots A_{\omega_{n-1}}\|_F
=|\mathbf{v}'  \mathbf{u}|^{\sum_{j=0}^r t_j} \Bigl( \prod_{j=1}^r 
\frac{ | \mathbf{v}' A_{\mathbf{w}_j} \mathbf{u}|} { |\mathbf{v}'  \mathbf{u}|}\Bigr)
\|A_{\bw_0}\mathbf{u}\|_2\cdot \|\mathbf{v}'A_{\bw_{r+1}^*}\|_2,
\end{equation}
Secondly, if $t_{r+1}^*>0$, we have
$$
 \|A_{\omega_0}\cdots A_{\omega_{n-1}}\|_F=|\lambda|^{\sum_{j=0}^{r} (t_j-1)+(t_{r+1}^*-1)}\Bigl( \prod_{j=1}^{r+1}  | \mathbf{v}' A_{\mathbf{w}_j} \mathbf{u}| \Bigr)\cdot 
\|A_{\bw_0}\mathbf{u}\|_2 \cdot\|\mathbf{v}\|_2
$$
which is then written as follows
\begin{equation}\label{expression2}
\|A_{\omega_0}\cdots A_{\omega_{n-1}}\|_F=|\mathbf{v}'  \mathbf{u}|^{\sum_{j=0}^{r} t_j+t_{r}^*} \Bigl( \prod_{j=1}^{r+1}
\frac{ | \mathbf{v}' A_{\mathbf{w}_j} \mathbf{u}|} { |\mathbf{v}'  \mathbf{u}|}\Bigr)
\|A_{\bw_0}\mathbf{u}\|_2\cdot \|\mathbf{v}\|_2.
\end{equation}

Having the expressions \eqref{expression1} and \eqref{expression2}, it is evident that
conditions (a) and (b) of Theorem \ref{thm:rank1} ensure that $\|A_{\omega_0}\cdots A_{\omega_{n-1}}\|_F\ne 0$ for all $n$ and that
$$
\lim_{n} \frac 1n \log \|A_{\bw_0}\mathbf{u}\|_2= \lim_{n} \log\frac 1n\|\mathbf{v}\|_2=  \lim_{n} \frac 1n \log  \|\mathbf{v}'A_{\bw_{r+1}^*}\|_2=0.
$$
So, what remain to prove are the following equalities
\begin{equation}\label{eq:freq}
\lim_{n\to\infty}\frac {\sum_{j=0}^{r} t_j+t_{r+1}^*}{n} =\rho_0
\end{equation}
\begin{equation}\label{eq:Mlim}
\lim_{n\to \infty}\frac 1{n}\sum_{j=1}^{r} \log \frac{ | \mathbf{v}' A_{\mathbf{w}_j} \mathbf{u}|} { |\mathbf{v}'  \mathbf{u}|}
= 
\sum_{\mathbf{w}\in \mathcal{R}}  F_{\mathbf{w}} \log \frac{|\mathbf{v}' A_{\mathbf{w}} \mathbf{u}|}
   {|\mathbf{v}'\mathbf{u} |}.
\end{equation}

The equality \eqref{eq:freq} holds just because
$ t_0+t_1 +\cdots +t_r +t_{r+1}^*$ is precisely the number of times the symbol $0$ appears
in $\omega|_0^{n-1}$.

 The equality \eqref{eq:Mlim} means
\begin{equation}\label{eq:Mlim2}
    \lim_{n\to \infty} \frac{1}{n} \sum_{j:\mathbf{w}_j \prec \omega|_0^{n-1}}  \log \left|\frac{\mathbf{v}' A_{\mathbf{w}_j} \mathbf{u}}{ \mathbf{v}'\mathbf{u}}\right| = 
       \sum_{\mathbf{w}\in \mathcal{R}}
F_\mathbf{w}  \log \left|\frac{\mathbf{v}' A_{\mathbf{w}} \mathbf{u}}{ \mathbf{v}'\mathbf{u}}\right| .
\end{equation}
Let us prove this limit.
For simplicity, for any $n\ge 1$ we denote 
$$
      S_n = \sum_{j:\mathbf{w}_j \prec \omega|_0^{n-1}}  \log \left|\frac{\mathbf{v}' A_{\mathbf{w}_j} \mathbf{u}}{ \mathbf{v}'\mathbf{u}}\right|, 
 $$
 For any integer $M\ge 1$, we decompose 
 \begin{equation}\label{eq:Sn}
 S_n=S_n'(M) +S_n''(M)
 \end{equation} with
 \begin{equation*}
  S_n'(M) = \sum_{j:\mathbf{w}_j \prec \omega|_0^{n-1}, \mathbf{w}_j \in \mathcal{R}_M}  \log \left|\frac{\mathbf{v}' A_{\mathbf{w}_j} \mathbf{u}}{ \mathbf{v}'\mathbf{u}}\right|,
  \end{equation*}
  \begin{equation*}
  S_n''(M) = \sum_{j:\mathbf{w}_j \prec \omega|_0^{n-1}, \mathbf{w}_j \in \mathcal{R}_M^c}  \log \left|\frac{\mathbf{v}' A_{\mathbf{w}_j} \mathbf{u}}{ \mathbf{v}'\mathbf{u}}\right|.
 \end{equation*}
 Since $\mathcal{R}_M$ is finite
 and the exact frequencies $F_{\mathbf{w}}$ exist, it is easy to see that
 $$
        \lim_{n\to\infty} \frac{S_n'(M)}{n} =\sum_{\mathbf{w}\in \mathcal{R}_M} F_\mathbf{w}\log \left|\frac{\mathbf{v}' A_{\mathbf{w}} \mathbf{u}}{ \mathbf{v}'\mathbf{u}}\right|.
 $$
 Then, from the decomposition \eqref{eq:Sn} we get
 $$
     \varliminf_{n\to \infty}  \frac{S_n}{n} =   \sum_{\mathbf{w}\in \mathcal{R}_M} F_\mathbf{w}  \log \left|\frac{\mathbf{v}' A_{\mathbf{w}} \mathbf{u}}{ \mathbf{v}'\mathbf{u}}\right|
     +  \varliminf_{n\to \infty}  \frac{S_n''(M)}{n},
 $$
 $$
     \varlimsup_{n\to \infty}  \frac{S_n}{n} =   \sum_{\mathbf{w}\in \mathcal{R}_M} F_\mathbf{w}  \log \left|\frac{\mathbf{v}' A_{\mathbf{w}} \mathbf{u}}{ \mathbf{v}'\mathbf{u}}\right|
     +  \varlimsup_{n\to \infty}  \frac{S_n''(M)}{n}.
 $$
 As 
 $$
     |S_n''(M)| \le  \sum_{j:\mathbf{w}_j \prec \omega|_1^n, \mathbf{w}_j \in \mathcal{R}_M^c}  \left|\log \left|\frac{\mathbf{v}' A_{\mathbf{w}_j} \mathbf{u}}{ \mathbf{v}'\mathbf{u}}\right|\right|,
 $$
 by the condition (c), we have
 $$
      \lim_{M\to \infty}   \varliminf_{n\to \infty}  \frac{S_n''(M)}{n} =  \lim_{M\to \infty}   \varlimsup_{n\to \infty}  \frac{S_n''(M)}{n}= 0.
 $$
 Thus we get 
 $$
     \varliminf_{n\to \infty}  \frac{S_n}{n} =    \varlimsup_{n\to \infty}  \frac{S_n}{n}= \lim_{M\to \infty}   \sum_{\mathbf{w}\in \mathcal{R}_M} F_\mathbf{w}  \log \left|\frac{\mathbf{v}' A_{\mathbf{w}} \mathbf{u}}{ \mathbf{v}'\mathbf{u}}\right|
     = \sum_{\mathbf{w}\in \mathcal{R}} F_\mathbf{w}  \log \left|\frac{\mathbf{v}' A_{\mathbf{w}} \mathbf{u}}{ \mathbf{v}'\mathbf{u}}\right|,
 $$
 where for the last equality we have used Lemma \ref{lem:lem2}. This finishes the proof of Theorem \ref{thm:rank1}.

\section{Product of non-negative matrices: Proof of Theorem \ref{thm:positive}} \label{Sect: Proof1}

In order to prove Theorem \ref{thm:positive} as a corollary of Theorem \ref{thm:rank1}, we need an investigation of  the structure of a generic sequence by studying its return words.

\subsection{Structure of generic sequence}

The proof of the following lemma, which  shows that the long return words are relatively sparse, is inspired by \cite{FW2022}.

\begin{lem}\label{lemma-long-small}
If $\omega$ is a $\nu$-generic point such that $\nu([0])>0$, 
we have 
$$
    \lim_{M\to \infty} \varlimsup_{n\to \infty} \frac{1}{n}\sum_{\mathbf{w}_j \prec \omega|_0^{n-1} |\mathbf{w}_j|>M}|\mathbf{w}_j|=0.
$$
\end{lem}

Let $U=\bigcup_{\mathbf{w}\in \mathcal{R}} [0\mathbf{w}0]$. Define inductively the return times to $U$ by
$$
   \tau_1 = \inf \{j\ge 0: \sigma^j(\omega) \in U\}; \quad  \tau_j = \inf \{j> \tau_{j-1}: \sigma^j(\omega) \in U\} \ \ (j \ge 2).
$$ 
In other words,
$$
       \forall j\ge 1, \quad \sigma^{\tau_j}(\omega) = 0{\mathbf{w}_j}0^{t_j} {\mathbf{w}_{j+1}}0^{t_{j+1}} \cdots.
$$
where $\mathbf{w}_j$ and ${t_j}$ are defined in \eqref{eq:omega}.   For $M\ge 1$, let 
$$
 \mathcal{R}_M= \{\mathbf{w}\in \mathcal{R}: |\mathbf{w}| >M\}.
$$

Consider the average length of long return words defined by
\[
S^M_i=\frac{1}{\tau_i}\sum_{j=1}^{i} \left |\mathbf{w}_j\right|{\bf 1}_{\mathcal{R}_M}(\mathbf{w}_j).
\]

What we have to prove is 
\begin{equation}\label{eq:LL=0}
\lim_{M\to \infty} \varlimsup_{i \to \infty} S_i^M =0.
\end{equation}
For any fixed $i\ge 1$, the average $S_i^M$  is decreasing in $M$,  and so is  $\varlimsup_{i\to\infty} S_i^M$. We shall prove \eqref{eq:LL=0} by contradiction.  Then suppose that there exists a $\delta>0$ such that 
$$
 \forall M\ge 1, \ \ \limsup_{i\to\infty} S_i^M\ge \delta.
 $$
 Let us fix an $M\ge 1$.  Let $\{i_\ell\}$ be a subsequence of
 integers  such that 
 \begin{equation}\label{eq:lemma-long-small 1} 
\lim_{\ell\to \infty} S_{i_\ell}^M\ge \delta.
 \end{equation}
 
For any $M\ge 1$ and $i\ge 1$ fixed,  consider the ``orbit measures along with long return words":
$$\nu_M^i:
       =\frac{1}{\tau_i}\sum_{j=1}^{i} {\bf 1}_{\mathcal{R}_M }(\mathbf{w}_j)\left(\sum_{n=1}^{|\mathbf{w}_j|-1}\delta_{\{ \sigma^n(0\mathbf{w}_j 0^{t_j}\mathbf{w}_{j+1}\cdots)\}}\right)
       $$
and the ``complementary orbit  measures" $\eta_M^i$ defined by 
$$\frac{1}{\tau_i}\sum_{n=0}^{\tau_i-1}\delta_{\sigma^n\omega} = \nu_M^i + \eta_M^i.$$
Recall that $\sigma^{\tau_j}\omega=0\mathbf{w}_j 0^{t_j}\mathbf{w}_{j+1}\cdots \in [0]$, but
$\sigma^n(0\mathbf{w}_j 0^{t_j}\mathbf{w}_{j+1}\cdots) \not\in [0]$
for all $1\le n\le |\mathbf{w}_j-1| $.
Thus the measure $\nu_M^i$ doesn't charge the cylinder $[0]$, i.e., 
\begin{equation}\label{eq:lemma-long-small 2}
\forall M\ge 1, \forall i, \quad 
\nu_M^i([0])=0.
\end{equation}
On the other hand, letting $Y=\overline{\{\sigma^n\omega\}_{n\ge 1}}$ be the orbit closure of $\omega$.
In view of  \eqref{eq:lemma-long-small 1} and the definition of $\nu_M^i$, we have
\begin{equation}\label{eq:lemma-long-small 3}
                      \forall M\ge 1, \quad
                      \nu_M^{i_\ell}(Y) \ge  \frac{1}{\tau_i}\sum_{j=1}^{i} {\bf 1}_{\mathcal{R}_M }(\mathbf{w}_j) (|\mathbf{w}_j|-1)\ge \frac{\delta}{2}\ \  \textrm{ if } \ell \gg 1.
\end{equation}
Up to taking a subsequence of $\{i_\ell\}$, we can assume that 
$$\nu_M^{i_\ell}\rightharpoonup \nu_M^\infty \ \textrm{ and }\ \eta_M^{i_\ell} \rightharpoonup \eta_M^\infty$$
where $\nu_M^\infty$ and $\eta_M^\infty$ are some measures concentrated on $Y$. Note that $\nu_M^\infty$ and $\eta_M^\infty$ are not necessarily probability measures.  
Because of \eqref{eq:lemma-long-small 2}, we  have 
\begin{equation}\label{eq:lemma-long-small 4}
\forall M\ge 1, \quad \nu_M^\infty([0])=0.
\end{equation}
By \eqref{eq:lemma-long-small 3}, we have 
\begin{equation}\label{eq:lemma-long-small 4-1}
\forall M\ge 1, \quad 
\nu_M^{\infty}(Y)\ge \frac{\delta}{2}.
\end{equation}
Since $\omega$ is $\nu$-generic, the measure $\nu_M^i+\eta_M^i$ converges to $\nu$ in the weak star topology, as $i\to\infty$. Thus we have 
\begin{equation}\label{eq:lemma-long-small 4-2}
\forall M \ge 1, \quad \nu_M^\infty+\eta_M^\infty=\nu.
\end{equation}

In the following, we will study in more details the measure  $\nu_M^\infty$.  Associated to the $j$-th return word $\mathbf{w}_j$, we consider the orbit probability measure
$$
     \pi_j := \frac{1}{|\mathbf{w}_j|-1}\sum_{n=1}^{|\mathbf{w}_j|-1} \delta_{\{ \sigma^n(0\mathbf{w}_j 0^{t_j}\mathbf{w}_{j+1}\cdots)\}}.
 $$
This is the orbit measure along with $\mathbf{w}_j$.  
Let 
$$A_M=\left\{\pi_j: \mathbf{w}_j\in \mathcal{R}_M\right\},$$
the set of  orbit measures corresponding to  long return words. Notice that
$$
    | \mathbf{w}_j | -1 \ge M \quad {\rm when}\ \ \pi_j \in A_M.
$$  

\begin{lem} \label{lem:lim-mu}
Suppose that we have a sequence of measures $\mu_t \in A_{M_{t}}$, with $M_t \to \infty$ as $t\to \infty$, such that $\mu_t  \rightharpoonup \mu$. Then $\mu$ is $\sigma$-invariant. 
\end{lem}
Indeed, it is a direct consequence of the following known fact:  Given two sequences of integers
$\{p_k\}, \{q_k\}\subset\mathbb{N}$  with $q_k-p_k\to\infty $ as $k\to\infty$ and a sequence of points $\{x_k\} \in \{0, 1, \cdots, m-1\}^\mathbb{N}$. If 
$$\frac{1}{q_k-p_k}\sum_{j=p_k}^{q_k}\delta_{\sigma^j(x_k)}\rightharpoonup \lambda, \textrm{ as } k\to\infty,$$   
then the limit measure $\lambda$ is a $\sigma$-invariant measure.
\medskip

Continue our discussion. Let $\mathcal{A}$ be the weak-* closure of $\bigcup_M A_M$,   a compact set of measures. Writing $\nu_M^i$ as 
$$\nu_M^i=\frac{1}{\tau_i}\sum_{j=1}^{i}  {\bf 1}_{\mathcal{R}_M}(\mathbf{w}_j) (|\mathbf{w}_{j}|-1) \pi_j ,$$
or equivalently 
$$\nu_M^i=\int_ {\mathcal{A}}\mu dQ_M^i(\mu)$$ 
where $Q_M^i$ is a discrete measure on $A_M (\subset \mathcal{A})$ with total mass not exceeding $1$.
Since $\nu_M^\infty$ is the weak-* limit of $\nu_M^{i_\ell}$,  
it can be written as 
\begin{equation}\label{eq:nu_M_infty}
\nu_M^\infty=\int_\mathcal{A} \mu dQ_M(\mu).
\end{equation}
for some  limit measure $Q_M$ of $Q_M^{i_\ell}$ (not necessarily probability measure) on the space $\overline{A}_M$. 
Let us explain the obtention of \eqref{eq:nu_M_infty}.
Recall that equipped with the weak-* topology  the dual space $C^*(\Omega)$ is  locally compact and metrizable.  The set $\mathcal{A}$ is a compact subset in $C^*(\Omega)$. 
That $Q_M^{i_\ell}$ converges to $Q_M$ means 
$$
     \int_{\mathcal{A}} \varphi(\mu) d Q_M^{i_\ell}(\mu) \to \int_{\mathcal{A}} \varphi(\mu) d Q_M(\mu)
$$ 
for all weak-$*$ continuous function $\varphi$. In particular, as the function 
$\mu \to \int_\Omega fd\mu$ is weak-$*$ continuous, we have 
$$
    \int_{\mathcal{A}} \int_\Omega f d\mu d Q_M^i(\mu) \to \int_{\mathcal{A}} \int_\Omega f d\mu d Q_M(\mu). 
$$
So, the equality (\ref{eq:nu_M_infty}) reads as 
$$
\int_\Omega f d\nu_M^\infty =  \int_{\mathcal{A}} \int_\Omega f d\mu d Q_M(\mu)
$$ 
for all continuous function $f$ on $\Omega$.


Take a sequence of integers $\{M_t\}$ tending to $\infty$ such that
$$\nu_{M_t}^\infty\rightharpoonup\nu_\infty^\infty, \quad
\eta_{M_t}^\infty\rightharpoonup\eta_\infty^\infty , \quad Q_{M_t}\rightharpoonup Q_\infty \quad \textrm{ as } t\to \infty.$$
Then from \eqref{eq:nu_M_infty}, we  have 
\begin{equation}\label{eq:nu-inf}
\nu_\infty^\infty=\int_\mathcal{A} \mu d Q_\infty (\mu);
\end{equation}
 from \eqref{eq:lemma-long-small 4}, we get 
\begin{equation}\label{eq:lemma-long-small 5}
\nu_\infty^\infty([0])=0;
\end{equation}
from \eqref{eq:lemma-long-small 4-2}, we get
\begin{equation}\label{eq:lemma-long-small 6}
\nu_\infty^\infty+\eta_\infty^\infty=\nu.
\end{equation}

Since $Q_\infty $ is the weak limit of $Q_{M_t}$,  it holds that for each $\mu\in {\rm supp}(Q_\infty)$ and each $r>0$, we have $Q_{M_t}(B(\mu,r))>0$ for all large enough $t$, where 
$B(\mu, r)$ is the ball centered at $\mu$ of radius $r$.  It follows that for every $\mu\in {\rm supp}(Q_\infty)$, there is a sequence of probability measures $\{\mu_t\}$ with $\mu_t\in A_{M_t}$ such that   $\mu_t\rightharpoonup \mu$, as $t\to\infty$. 
By Lemma \ref{lem:lim-mu},
the measure $\mu$ is $\sigma$-invariant on $Y$.   Then $\nu_\infty^\infty$ is $\sigma$-invariant, by (\ref{eq:nu-inf}).

By \eqref{eq:lemma-long-small 6}, the invariant measure $\nu_\infty^\infty$ is absolutely continuous with respect to the measure $\nu$ which is assumed ergodic. It follows that
$\nu_\infty^\infty= c \nu$ for some constant $c$, which is not zero by (\ref{eq:lemma-long-small 4-1}).
Finally,  (\ref{eq:lemma-long-small 5}) contradicts $\nu([0])>0$.

\subsection{Proof of Theorem \ref{thm:positive}}
If $\mathbf{v} A_{\mathbf{w}_j} \mathbf{u}'=0$ for some $j$, we have $L(\omega)=-\infty$ by Theorem \ref{thm:rank1}.
Assume that $\mathbf{v} A_{\mathbf{w}_j} \mathbf{u}'\not=0$ for all $j$. We have only to check the conditions (b) and (c) in Theorem \ref{thm:rank1}.

First remark that if $\Delta$ denote the maximum of the operator norms
$\|A_1\|, \cdots, \|A_{m-1}\|$ we have the  upper bound
$$
   |\mathbf{v} A_{\mathbf{w}_j} \mathbf{u}'| \le  |\mathbf{v}\| \|\mathbf{u}'\| \Delta^{|\mathbf{w}_j|}.
$$
On the other hand, from the definition of matrix product,  we  have the lower bound
$$
    |\mathbf{v} A_{\mathbf{w}_j} \mathbf{u}'| \ge a_* b_* \delta^{|\mathbf{w}_j|}
$$
where  $a_*$ (resp. $b_*$) is the minimum of the non zero entries of $\mathbf{v}$ (resp. $\mathbf{v}$) and $\delta$ is the minimum of all non zero entries of 
$A_1, \cdots, A_{m-1}$.
From the above two estimates we get 
$$
     |\log   |\mathbf{v} A_{\mathbf{w}_j} \mathbf{u}'| | \le C |\mathbf{w}_j|
$$
for some constant $C>0$. It follows that 
$$
   \frac{1}{n}\sum_{\mathbf{w}_j \prec \omega|_0^{n-1}, |\mathbf{w}_j|>M} |\log |\mathbf{v}'A_{\mathbf{w}_j} \mathbf{u}|| 
   \le   \frac{C}{n}\sum_{\mathbf{w}_j \prec \omega|_0^{n-1}, |\mathbf{w}_j|>M} |\mathbf{w}_j|. 
$$
But by Lemma \ref{lemma-long-small}, the right hand term in the above inequality tends to zero as $n \to \infty$. Thus we have checked the condition
(c) in Theorem \ref{thm:rank1}.

Similarly, we can deduce that $\gamma_n = O(|\mathbf{w}_{r(n) +1}|)$. Again, we get $|\mathbf{w}_{r(n) +1}| =o(n)$ by Lemma \ref{lemma-long-small}.
Thus we have checked the condition
(b) in Theorem \ref{thm:rank1}.

\section{Product of random matrices: Proof of Theorem \ref{thm:random}} \label{Sect: Proof2}

Let us first identify the return words for a $\mu$-typical sequence $\omega$. Notice that if 
$\mu([0\mathbf{w}0])=0$, then
$$
    \mu\left(\bigcup_{k=0}^\infty \sigma^{-k}([0\mathbf{w}0])\right) \le \sum_{k=0}^\infty \mu(\sigma^{-k}([0\mathbf{w}0]))=0.
$$
So, for almost every $\omega$, the word $[0\mathbf{w}0]$ never appears in $\omega$. On the other hand,  every  word $0\mathbf{w}'0$ such that 
$\mu(0\mathbf{w}'0)>0$ appears infinitely many times  in $\omega$ by Poincar\'e's recurrence theorem. Thus  the set $\mathcal{R}$ stated in Theorem \ref{thm:random}
is the set of return words for almost every $\omega$.  
Hence, such $\omega$'s have the decomposition \eqref{eq:omega}. 
We will consider such typical points $\omega$.
\medskip

{\em Case I.}
\underline{Assume $\mathbf v'A_\mathbf{w}\mathbf u=0$ for some $\mathbf{w}\in \mathcal{R}$}.  We immediately conclude that $L=-\infty$ by Theorem \ref{thm:rank1}.
\medskip 

{\em Case II.} \underline{From now on we assume that $\mathbf v'A_\mathbf{w}\mathbf u\not=0$ for all $\mathbf{w}\in \mathcal{R}$}.
We are going to apply Theorem \ref{thm:rank1}, by checking that all the conditions (a)-(c) in Theorem \ref{thm:rank1} are satisfied for
$\mu$-almost all  $\omega$.  
\medskip

{\it The condition (a) is almost surely satisfied.} 
It suffices to check that $A_{\mathbf{w}_0} \mathbf u\cdot \mathbf v' A_{\mathbf{w}^*_{t_{r+1}}} \not=0$. Indeed, otherwise $A_0=\mathbf u\cdot \mathbf v'=0$, 
because both matrices $A_{\mathbf{w}_0}$ and $ A_{\mathbf{w}^*_{t_{r+1}}}$ are invertible. 
\medskip

{\it The condition (b) is almost surely satisfied.}   
It is essentially to show that with the notation used in 
\eqref{eq:omega-n}, we have 
\begin{equation}\label{eq:kac1}
a.s. \quad t_r+|w_r|=o(n).
\end{equation}
To this end, let us consider the induced dynamics.
Consider the disjoint union over return words $\mathbf{w}\in \mathcal{R}$ of cylinders $[0\mathbf{w}0]$:
      $$U = \bigcup_{\mathbf{w} \in \mathcal{R}}[0\mathbf{w}0].$$
Consider the first return time to U defined by
     $$\tau(x) = \inf \{n\ge 1: \sigma^n x \in U \}.$$
It is the first time to see a return word $\mathbf{w}$ (or more exactly $0\mathbf{w}0$). 
By the ergodicity of $\mu$ and the fact that $\mu(U)>0$, we get $\tau(x)<\infty$ for almost all points $x$. 
Then let us consider the induced map $T: U \to U$ defined by 
 $$T(x) = \sigma^{\tau(x)}(x).$$
 The restricted measure $\mu|_U$ is still $T$-invariant and ergodic. 
Define higher order return times, inductively,
by
    $$\tau_1(x) =\tau(x), \qquad 
    \tau_{k} (x) = \tau(T^{k-1}(x)).$$
Notice that  $\tau_k(\omega)$ is the time difference when we see the word $0\mathbf{w}_{k-1}0$ and the the word $0\mathbf{w}_{k}0$ , where $\mathbf{w}_j'$s denotes  the $j$-th return word in the decomposition of $\omega$, see \eqref{eq:omega}. So, we have
\begin{equation}\label{eq:kac2}
   a.s. \quad   t_k+|w_k| =O(\tau_k(\omega)),
\end{equation}
which implies $|w_k| =O(\tau_k(\omega))$ holds almost surely. This fact will be useful. 
By Kac's theorem, $\tau$ is integrable and even $\int_U\tau(x) d\mu=1$. Then, by the Birkhoff theorem  applied to the induced system, the Birkhoff average
     $$\frac{1}{k} \sum_{j=1}^k \tau_j(x) =  \frac{1}{k} \sum_{j=0}^{k-1}  \tau(T^jx)
     $$
converges a.e. Consequently, $ \tau(T^{k-1}x)=o(k)$ a.s. 
This, together with \eqref{eq:kac2},  implies  \eqref{eq:kac1}.       Now, let us check the condition (b). From the estimates
$$
   \|\mathbf{v}'A_{\mathbf{w}'}\| \le \|\mathbf{v} \|\|A_{\mathbf{w}'}\|, \qquad
    \|\mathbf{v}\| =
     \|\mathbf{v}'A_{\mathbf{w}'} A_{\mathbf{w}'}^{-1}\|
     \le  \|\mathbf{v}'A_{\mathbf{w}'}\|\| A_{\mathbf{w}'}^{-1}\|
$$
we get 
$$
\frac{\|\mathbf{v}\|}{\|A_{\mathbf{w}'}^{-1}\|}\le \|\mathbf{v}'A_{\mathbf{w}'}\| \le \|\mathbf{v} \|\|A_{\mathbf{w}'}\|
$$
But both $\|A_{\mathbf{w}'}\|$ and $\|A_{\mathbf{w}'}^{-1}\|$
increase at most exponentially with $|\mathbf{w}|$. So, the quantity $\gamma_n$
involved in the condition (b) satisfies
$$
     \gamma_n =O(|\mathbf{w}_{r+1}|).
$$
This implies $\gamma_n  = o(n)$ a.s, with the help of \eqref{eq:kac1}. Thus we have checked the condition (b) for almost all
$\omega$.
\medskip

   {\it The condition (c) is almost surely satisfied.}  
   We do not check the condition (c). 
As we have shown in the proof of Theorem \ref{thm:rank1},
what we have to prove is that \eqref{eq:Mlim2} holds almost surely, namely 
\begin{equation}\label{eq:Mlim3}
  a.s. \quad   \lim_{n\to \infty} \frac{1}{n} \sum_{j:\mathbf{w}_j \prec \omega|_0^{n-1}}  \log \left|\frac{\mathbf{v}' A_{\mathbf{w}_j} \mathbf{u}}{ \mathbf{v}'\mathbf{u}}\right| = 
       \sum_{\mathbf{w}\in \mathcal{R}}
F_\mathbf{w}  \log \left|\frac{\mathbf{v}' A_{\mathbf{w}} \mathbf{u}}{ \mathbf{v}'\mathbf{u}}\right| .
\end{equation}
Notice that $$ 
\frac{1}{n} \sum_{j:\mathbf{w}_j \prec \omega|_0^{n-1}}  \log \left|\frac{\mathbf{v}' A_{\mathbf{w}_j} \mathbf{u}}{ \mathbf{v}'\mathbf{u}}\right| 
=
\frac{1}{n}\sum_{k=0}^{n-1} h(\sigma^k \omega),
$$
where 
$$
   h(\omega) = \sum_{\mathbf{w}\in \mathcal{R}} 1_{[0\mathbf{w}0]} (\omega) \log \left|\frac{\mathbf{v}' A_{\mathbf{w}} \mathbf{u}}{ \mathbf{v}'\mathbf{u}}\right|. 
$$
Notice that the right hand side of \eqref{eq:Mlim3}
is nothing by the mean value of $h$. That is to say
$$
     \int h(\omega) d\mu(\omega) = \sum_{\mathbf{w}\in \mathcal{R}}
F_\mathbf{w}  \log \left|\frac{\mathbf{v}' A_{\mathbf{w}} \mathbf{u}}{ \mathbf{v}'\mathbf{u}}\right|.
$$
So, \eqref{eq:Mlim3} is a consequence of the Birkhoff ergodic theorem. To finish the proof, we point out that $h$
is not necessarily integrable. But the positive part of $h$ is integrable, because of Lemma \ref{lem:lem1} and  the fact $$\log |\mathbf{v}'A_{\mathbf{w}}\mathbf{u}| \le C |\mathbf{w}|$$
for some constant $C>0$.
Therefore the ergodic theorem is applicable in this case.

  \section{Computation of exact frequencies} \label{Sect: ExactF}
  
It was pointed out  that the exact frequency $F_\mathbf{w}$, which appears in the our main formula \eqref{eq:LiapRank1},  is  nothing but the the
   asymptotic frequency of the word $0\bw 0$, namely 
   $F_\mathbf{w}=\rho_{0\mathbf{w}0}$ (cf. \eqref{eq:F=rho}).
   But it is possible to compute the exact frequency $F_\mathbf{w}$ by the asymptotic frequencies of return words
   $\mathbf{w}'$ with their lengths not exceeding that of $\mathbf{w}$. This remark is useful because the asymptotic frequencies of return words of shorter lengths are easier to compute. 
   
 \subsection{Computation of $F_\mathbf{w}$ by $\rho_{\mathbf{w}'}$ with $|\mathbf{w}'|\le |\mathbf{w}|$}

  Let us assume that $\omega$ is $\nu$-generic for some shift invariant measure $\nu$ and the collection of return words $\mathcal R$ in $\omega$ is finite. Here we refer to the decomposition \eqref{eq:decomp2} of $\omega$.
  We will compute the  exact frequencies $F_{\mathbf{w}}$ ($\mathbf{w}\in \mathcal{R}$) in terms of the cylinder measures $\nu([\mathbf{w'}])$ ($\mathbf{w}'\in \mathcal{R}$) and the numbers $N_{\mathbf{w}}(\mathbf{w}')$ 
  with $\mathbf{w}, \mathbf{w}' \in \mathcal{R}$, which are defined as follows.  
  
  For two given words $\mathbf{w}$ and $\mathbf{w}'$ with $\mathbf{w} \prec \mathbf{w}'$ (meaning that $\mathbf{w}$ is a proper subword of $\mathbf{w}'$), we denote by $N_{\mathbf{w}}(\mathbf{w}')$  the number of occurrences of $\mathbf{w}$ in $\mathbf{w}'$ when reading $\mathbf{w}'$ from left to right. For example, $N_{01}(101001101)=3$.
  In practice,  for a given finite collection of return words $\mathcal R$, it is easy to compute
   the numbers $N_{\mathbf{w}}(\mathbf{w}')$. 
   
  As $\mathcal{R}$ is assumed finite. Let $\ell$ be the maximal length of return words, i.e.
  $$
      \ell=\max\{|\mathbf{w}|: \mathbf{w} \in \mathcal{R}\}.
  $$
For any return word $\mathbf{w} \in \mathcal{R}$, put  \begin{eqnarray*}
S_{\mathbf{w}}^{(0)} & = & \nu([\mathbf{w}]), \\
     S_{\mathbf{w}}^{(j)} & = & \sum_{\mathbf{w} \prec \mathbf{w} '\prec \cdots \prec \mathbf{w}^{(j)}} N_{\mathbf{w} }(\mathbf{w} ')N_{\mathbf{w}' }(\mathbf{w} '')\cdots N_{\mathbf{w}^{(j-1)} }(\mathbf{w} ^{(j)})
      \nu([\mathbf{w} ^{(j)}]),  \ \  (1\le j \le \ell -|\mathbf{w}|)\\
     S^{(j)}_{\mathbf{w}}  & =& 0,  \ \ (j>\ell -|\mathbf{w}|).
  \end{eqnarray*}
  where the sum is taken over all $j$-uple return words $(\mathbf{w}', \mathbf{w}'',\cdots, \mathbf{w}^{(j)})$ such that $\mathbf{w}^{(i)}$ is a proper subword of $\mathbf{w}^{(i+1)}$
  for $0\le i <j$ (with convention $\mathbf{w}^{(0)}=\mathbf{w}$). 
 
 The numbers $S_{\mathbf{w}}^{(j)}$, mixing the numbers $N_{\mathbf{w}}(\mathbf{w}')$ and the cylinder measures $\nu([\mathbf{w}])$, allow us to compute the exact frequencies $F_{\mathbf{w}}$,
 as stated in the following theorem.
  
  \begin{thm}\label{thm:ExactF} Let $\omega$   is decomposed as in (\ref{eq:omega}). Suppose that the set $\mathcal{R}$ of return words is finite
  and the frequency of each $\mathbf{w}\in \mathcal{R}$ exists, which is denoted by $\nu([\mathbf{w}])$.
  Then for any return word $\mathbf{w}\in \mathcal{R}$,  the  exact frequency $F_{\mathbf{w}} $ exists and is equal to
  \begin{equation}\label{eq:Fw}
      F_{\mathbf{w}} = \sum_{j=0}^{\ell -|\mathbf{w}|} (-1)^j S^{(j)}_{\mathbf{w}} = \sum_{j\ge 0} (-1)^j S^{(j)}_{\mathbf{w}}.
        \end{equation}
  \end{thm}

 \subsection{Proof of Theorem \ref{thm:ExactF}} 
 
  We prove the formula (\ref{eq:Fw}) by induction on the length of $\mathbf{w}$ in decreasing order. 
 We start with the remark that for any return word $\mathbf{w} \in \mathcal{R}$, we have 
\begin{equation}\label{eq:3b}
    \nu ([\mathbf{w}]) = F_{\mathbf{w}} + \sum_{\mathbf{w}' \in \mathcal{R}, |\mathbf{w}'|>|\mathbf{w}|}  N_{\mathbf{w}}(\mathbf{w}') F_{\mathbf{w}'},
\end{equation}
under the condition that the exact frequencies $F_{\mathbf{w}'}$ exist for all return words $\mathbf{w}'$ containing $\mathbf{w}$ as proper subword. 
This is a consequence of the  existences of frequencices:
 \begin{equation}\label{eq:Freq_w}
  \forall \mathbf{w} \in \mathcal{R}, \quad \nu([\mathbf{w}]) = \lim_{n\to \infty} \frac{1}{n}\sum_{k=0}^{n-1}1_{[\mathbf{w}]}(\sigma^k \omega).
\end{equation}
Indeed,  when we read the prefix of  $\sigma^k\omega$ and see $\mathbf{w}$  (recall here the decomposition (\ref{eq:omega}) of $\omega$), either we meet the return word $\mathbf{w}$ in the decomposition
or we meet longer return words $\mathbf{w}'$ which contain $\mathbf{w}$ as subword. The summands in the  right hand side of (\ref{eq:Freq_w})  which correspond to the return words $\mathbf{w}$ that we meet  gives the contribution $F_{\mathbf{w}}$ in (\ref{eq:3b}). While when we meet $\mathbf{w}'$, we will read from left to right   a number $N_{\mathbf{w}}(\mathbf{w}')$ of $\mathbf{w}$'s in each $\mathbf{w}'$. 
Since the exactly frequency of $\mathbf{w}'$ is $F_\mathbf{w}'$,   we get the other contribution of the right hand side of (\ref{eq:Freq_w}), the second term on the right hand side of  (\ref{eq:3b}).
\medskip

Now let us start our proof by induction. 
 When $|\mathbf{w}|=\ell$, the formula (\ref{eq:3b}) reads  as $\nu([\mathbf{w}])=F_{\mathbf{w}}$, which  confirms the formula  (\ref{eq:Fw}) in this case.

Assume now $\mathbf{w} \in \mathcal{R}$ with $|\mathbf{w}|=\ell-1$. Then the formula (\ref{eq:3b}) reads as
   $$
    \nu ([\mathbf{w}]) = F_{\mathbf{w}} + \sum_{\mathbf{w}' \in \mathcal{R}, |\mathbf{w}'|>|\mathbf{w}|}  N_{\mathbf{w}}(\mathbf{w}') \nu([\mathbf{w}']).
    $$
    because $|\mathbf{w}'|= \ell$, and $F_{\mathbf{w}'}=\nu([\mathbf{w}'])$ as we have just seen. Thus
    \begin{equation}\label{eq:l-1}
       F_{\mathbf{w}} = \nu([\mathbf{w}]) - \sum_{\mathbf{w}'\in \mathcal{R}: |\mathbf{w}'|>|\mathbf{w}|} N_{\mathbf{w}}(\mathbf{w}')  \nu([\mathbf{w}']) = S_{\mathbf{w}}^{(0)} - S_{\mathbf{w}}^{(1)}.
\end{equation}
 So the formula  (\ref{eq:Fw}) is also confirmed  when  $|\mathbf{w}|=\ell-1$.
 
  Suppose  $F_{\mathbf{w}}$ exists and the formula  (\ref{eq:Fw}) is true when $|\mathbf{w}|\ge k+1$. For $\mathbf{w}\in \mathcal{R}$ with $|\mathbf{w}|= k$,
 by (\ref{eq:3b}), we get
 $$
    F_{\mathbf{w}} = \nu ([\mathbf{w}]) -  \sum_{\mathbf{w}' \in \mathcal{R}, |\mathbf{w}'|>|\mathbf{w}|}  N_{\mathbf{w}}(\mathbf{w}') F_{\mathbf{w}'}.
 $$
 By the induction hypothesis,   we have 
 $$
    F_{\mathbf{w}'} = \sum_{j\ge 0} (-1)^j S^{(j)}_{\mathbf{w}'}.
 $$
Thus, by substituting this into the preceding equality we get  
  \begin{eqnarray*}
   F_{\mathbf{w}} &=& \nu ([\mathbf{w}]) -  \sum_{ |\mathbf{w}'|>|\mathbf{w}|}  N_{\mathbf{w}}(\mathbf{w}')  \sum_{j\ge 0} (-1)^j S^{(j)}_{\mathbf{w}'}\\
       & = &  \nu ([\mathbf{w}]) +  \sum_{j\ge 0}    \sum _{ |\mathbf{w}'|>|\mathbf{w}|} (-1)^{j+1} N_{\mathbf{w}}(\mathbf{w}')  S^{(j)}_{\mathbf{w}'}.
   \end{eqnarray*}
   Recall that
   $$
         S^{(j)}_{\mathbf{w}'} =
        \sum_{\mathbf{w}'\prec \mathbf{w} ''\prec \cdots \prec \mathbf{w}^{(j+1)}} N_{\mathbf{w}' }(\mathbf{w} '')\cdots N_{\mathbf{w}^{(j)} }(\mathbf{w} ^{(j+1)})
      \nu([\mathbf{w} ^{(j+1)}])
  $$
  where the sum is taken over $\mathbf{w} '', \cdots, \mathbf{w}^{(j+1)}$ ($\mathbf{w} '$ being fixed). So,
 \begin{eqnarray*}
  F_{\mathbf{w}} &=& \nu ([\mathbf{w}]) + \sum_{j\ge 0} (-1)^{j+1}
       \hspace{-2em}\sum_{|\mathbf{w}|<|\mathbf{w}'|<|\mathbf{w} ''|<\cdots <|\mathbf{w}^{(j+1)}|}  N_{\mathbf{w}}(\mathbf{w}')  N_{\mathbf{w}' }(\mathbf{w} '')\cdots N_{\mathbf{w}^{(j)} }(\mathbf{w} ^{(j+1)})  \nu([\mathbf{w} ^{(j+1)}])\\
  &=& \nu ([\mathbf{w}]) + \sum_{k\ge 1} (-1)^k S^{(k)}_{\mathbf{w}}.
 \end{eqnarray*}
 Recall that  $\nu ([\mathbf{w}])= S^{(0)}_{\mathbf{w}}$.  We have thus finished the proof of Theorem \ref{thm:ExactF}. 
 \medskip

 \subsection{Examples}\, \ \ 

In the following, we consider some sequences $\omega$ of $0$ and $1$, which are generic for some invariant measure $\nu$.
\medskip

 {\bf Example 1.} Assume that $\mathcal{R}=\{1\}$. By definition, we have $F_1=\nu([1])$. It is the case for the Fibonacci sequence that will be examined later.
\medskip 
 
  {\bf Example 2.} Assume that $\mathcal{R}=\{1, 11\}$. Then
  $$
   S_{11}^{(0)}=\nu([11]), \ \ \ S_{11}^{(j)}=0 \ \ \text{for} \ j\ge 1
  $$
  and 
  $$
   S_{1}^{(0)}=\nu([1]), \ \ \ S_1^{(1)}= 2 \nu([11]),  \ \ \  S_{11}^{(j)}=0 \ \ \text{for} \ j\ge 2.
  $$
  Therefore the formula \eqref{eq:Fw} gives us
  \begin{equation}\label{eq:}
  F_{11}=\nu([11]), \quad F_1=\nu([1])-2 \nu([11])= \nu([10])-\nu[11].
  \end{equation}
  Here we have used the fact $N_1(11)=2$. This is the case
  for he Thue-Morse sequence, that will be examined later.
  \medskip

   {\bf Example 3.} Assume that $\mathcal{R} = \{1, 11, 111\}$. It is the case for the characteristic function of the square free integers, that will be examined later.  Then by Theorem \ref{thm:ExactF}, we have 
  \begin{equation}\label{ex:1}
       F_{111} = \nu([111]), \quad F_{11} = \nu([11]) - 2 \nu([111]), \quad  F_{1} = \nu([1]) - 2\nu([11]) + \nu([111]).
       \end{equation}
Indeed, the first equality is trivial. The formula \eqref{eq:Fw} implies immediately the second one:
$$
F_{11} =  \nu([11]) -\sum_{|11|<|\mathbf{w}'|} \nu([\mathbf{w}'])  = \nu([11])-N_{11}(111)\nu([111]) = \nu([11])- 2\nu([111])
$$
because $\mathbf{w}'=111$ is the unique choice over which we take the sum. Now we state the formula (\ref{eq:Fw}) for $F_1$:
\begin{eqnarray*}
   F_1 & = & \nu([1]) -\sum_{|1|<|\mathbf{w}'|} N_1(\mathbf{w}')\nu([\mathbf{w}']) + \sum_{|1|<|\mathbf{w}'|<|111|} N_1(\mathbf{w}') N_{\mathbf{w}'}(111)\nu([111]).
\end{eqnarray*}
In the first sum, $\mathbf{w}'$ may be $11$ with $N_1(11)=2$ or $111$ with $N_{1}(111)=3$; in the second sum,  
$\mathbf{w}'$ must be $11$ and $N_1(11)=N_{11}(111)=2$. So
$$
   F_1 = \nu([1]) - (2\nu([11])+3\nu([111])) + 4\nu([111]) =  \nu([1]) - 2\nu([11])+ \nu([111]). 
$$

\section{Matrices selected by substitutive sequences}\label{sect:SS}
 Given $m$ matrices $A_0, A_1, \cdots, A_{m-1}$ ($m\ge 2$)
 with $A_0$ being assumed of rank one.
 We select a sequence of matrices $(A_{\omega_n})_{n\ge 1}$
 by some sequence $\omega=(\omega_n) \in \{0,1, \cdots, m-1\}^\infty$. 
 Typical primitive substitutive sequences $\omega$ include the famous Thue-Morse sequence, Fibonacci sequence and Rudin-Shapiro sequence, etc.
  For a sequence of matrices $(A_{\omega_n})_{n\ge 1}$ selected by a primitive substitutive sequence $(\omega_n)$, we can  use Theorem
  \ref{thm:rank1} 
 to compute the Lyapunov exponent of $$\lim \frac{1}{n}\log\|A_{\omega_0}A_{\omega_1}\cdots A_{\omega_{n-1}}\|$$  
 if $A_1, \cdots, A_{m-1}$ are assumed non-negative 
 (cf. Theorem \ref{thm:positive}) or invertible (cf. Theorem \ref{thm:random}).  
 We assume that $0$ appears infinitely in $(\omega_n)$.
 We have to compute the frequencies of the return words of $0$. 

 In this section, we compute Lyapunov exponents when $\omega$ is a substitutive sequence. 
 Let us recall basic facts about substitutions and the method of computation of frequencies (cf. \cite{Queffelec}). 
 
 \subsection{Basic notions and notation}
 Let $\mathcal{A}=\{0,1, \cdots, s-1\}$ be an alphabet  of $s\ge 2$ letters.  Let $\mathcal{A}^*=\bigcup_{n=0}^\infty \mathcal{A}^n$ be the set of all words on the alphabet $\mathcal{A}$, with $\mathcal{A}^0=\{\emptyset\}$
 consisting of the empty word. The length $|w|$ of a word $w=w_0w_1\cdots w_{n-1}\in \mathcal{A}^n$ is defined to be $n$. Two words $u, v\in \mathcal{A}^*$ can be concatenated to a new 
 word $uv$. Equipped with the concatenation, $\mathcal{A}^*$ becomes a monoid with $\emptyset$ as identity.  
 
 Let $\mathcal{A}^+=\mathcal{A}^*\setminus\{\emptyset\}$. A {\em substitution} $\zeta$ on $\mathcal{A}$ is a map from $\mathcal{A}$ into $\mathcal{A}^+$. It induces a morphism of the monoid $\mathcal{A}^*$ as follows
 $$
     \zeta(\emptyset)=\emptyset; \quad \zeta(w) = \zeta(w_0)\zeta(w_1)\cdots \zeta(w_{n-1}) \ \ {\rm for}\ \ w=w_0w_1\cdots w_{n-1} \in \mathcal{A}^n \subset \mathcal{A}^+.
 $$
 The {\em language} of the substitution $\zeta$ is defined by 
 $$
    \mathscr{L}_\zeta=\{{\rm subwords\ of}  \ \zeta^n(\alpha): \alpha\in \mathcal{A}, n\ge 0 \}.
 $$ 
 The language  $\mathscr{L}(x)$ of an infinite word $x\in \mathcal{A}^\mathbb{N}$ is defined to be the set of finite subwords of $x$.  The {{\em subshift} of the substitution $\zeta$ is defined to be
 $$
  X_\zeta=\{x\in \mathcal{A}^\mathbb{N}: \mathscr{L}(x) \subset  \mathscr{L}_\zeta\}.
 $$
 The substitution also defines a continuous map $\zeta: \mathcal{A}^\mathbb{N} \to \mathcal{A}^\mathbb{N}$: 
 $$
    \zeta(x) = \zeta(x_0)\zeta(x_1)\cdots \quad {\rm for} \ \ x=x_0x_1\cdots \in \mathcal{A}^\mathbb{N}.
 $$
 We are interested in the possible fixed points of this continuous map.  Now we make the following assumptions:
 \begin{itemize}
 \item[(i)]  For every $\beta\in \mathcal{A}$, $\lim_{n\to \infty}|\zeta^n(\beta)|=+\infty$;
 \item[(ii)]  For some $\alpha\in \mathcal{A}$, $\zeta(\alpha)$ has $\alpha$ as its first letter. 
 \end{itemize}
 Under these assumptions, $\zeta$ admits a fixed point $u$, called {\em substitutive sequence},  which is equal to $\zeta^\infty(\alpha)$, namely 
 $\lim_{n\to \infty} \zeta^n(\alpha)$. We also assume that every letter in $\mathcal{A}$ appears in the substitutive sequence $u$. The language of $u$ is simply
 $$
    \mathscr{L}(u) = \{{\rm subword\ of\ } \zeta^n(\alpha): n \ge 0 \}.
 $$ 
 A {\em morphic sequence} is the image of a substitutive sequence under a letter-to-letter projection.  
 
 Given two words $u$ and $v$, we denote by $L_u(v)$ the number of occurrences of $u$ in $v$. 
 For a substitution $\zeta$, its $\zeta$-{\em matrix} or {\em composition matrix}, is defined to be
 $$M_\zeta:=M:=(m_{\alpha, \beta})\ \ \  {\rm with }  \ \ \ m_{\alpha, \beta}=L_\alpha(\zeta(\beta)).$$  
 It is clear that $$\sum_\alpha m_{\alpha, \beta} = |\zeta(\beta)|.$$
  It is also easy to see that 
 $$M_{\zeta^n}=M_\zeta^n.$$
 
 The $\zeta$-{\em matrix} 
 doesn't contain all informations about the substitution, but some important informations. 
 First we introduce the {\rm composition function} $L:\mathscr{L}_\zeta \to \mathbb{R}^s$, defined by $$
 L(w) := (L_\alpha(w))_{\alpha \in \mathcal{A}},
 $$ 
 which counts the number of each letter contained in a given word $w$. The composition matrix  involves in the following linear recursive relation
 (which can be called {\em composition equation})
 $$
     L(\zeta(w)) = M\cdot L(w),
 $$
 which describes the evolution of the numbers of letters contained in  a word under the substitution.
 It follows that $L(\zeta^n(w)) = M^n\cdot L(w)$. This relation allows us to study the frequencies of letters contained in a substitutive sequence, and similar properties. 
 The Perron-Frobenius theorem can provide a complete solution. This is one of dynamical aspects of the substitution.
 
 Given a substitution $\zeta$ with substitutive sequence $u$. There are two associated subsystems of the shift $T: \mathcal{A}^\mathbb{N}\to \mathcal{A}^\mathbb{N}$. 
 One is $(X_\zeta, T)$ and the other is $(\overline{O(u)}, T)$, where $O(u)=\{T^n u \}_{n\ge 0}$ is the orbit of $u$. 
 Clearly $\overline{O(u)} \subset X_\zeta$.  But  when $\zeta$ is primitive, we have   $\overline{O(u)} = X_\zeta$. In the following, we make the following assumption
 \begin{itemize}
 \item[(iii)] $\zeta$ is primitive, namely $M^N >0$ for some $N\ge 1$.
 \end{itemize}
 If $\zeta$ is primitive, we even have $X_{\zeta^n}=X_\zeta$ for all $n\ge 1$ and the system $(X_\zeta, T)$ is minimal and uniquely ergodic \cite{Michel1974} (see also \cite{Queffelec}, p. 141).
 \subsection{The unique ergodic measure}
 Let $\nu$ be the unique invariant probability measure of the primitive substitutive system $(X_\zeta, T)$. 
  We are  interested in the computation of the 
 frequencies of patterns in the  substitutive sequence $u$ (cf. Section 5.4 of \cite{Queffelec}). 
 
 Let $a_0a_1\cdots a_{t-1}$ be a given pattern, namely a finite word. Recall that $[a_0a_1\cdots a_{t-1}]$ denotes the cylinder set consisting of all
 infinite words having $a_0a_1\cdots a_{t-1}$ as prefix.  The frequency of  $a_0a_1\cdots a_{t-1}$ is nothing but
 the measure $\nu( [a_0a_1\cdots a_{t-1}])$, according to the unique ergodicity. That is to say
 $$
     \lim_{n\to \infty}\frac{1}{n}\sum_{k=0}^{n-1} 1_{[a_0a_1\cdots a_{t-1}]}(T^k u) = \nu([a_0a_1\cdots a_{t-1}]).
 $$ 
 When $t=1$,  if $\rho$ denotes the maximal eigenvalue  of $M$ (called the Perron-Frobenius eigenvalue), from the composition equation $L(\zeta^n(w)) = M^n\cdot L(w)$
 and the Perron-Frobenius theorem, we get that $\rho^{-n} L(\zeta^n(\alpha))$ is asymptotic to a Perron-Frobenius vector. Normalizing this eigenvector to a probability vector,
 we get the frequencies of letters contained in $u$. 
 
 The frequencies of patterns of length $t>1$ can be similarly computed, through an induced substitution. 
 Let $\mathcal{A}_t$ be the set of all words of length $t$ contained in $u$. We consider the induced substitution $\zeta_t$ on the alphabet $\mathcal{A}_t$ defined as follows. 
 For $\omega=\omega_0\omega_1\cdots \omega_{t-1}\in \mathcal{A}_t$, assume that 
 $$
\zeta(\omega)= x_0x_1\cdots x_{|\zeta(\omega_0)|-1}  x_{|\zeta(\omega_0)|} \cdots x_{|\zeta(\omega_0)|+|\zeta(\omega_1)|+1} \cdots x_{|\zeta(\omega)|-1}.
 $$  
 Then we define
 $$
   \zeta_t(\omega) = (x_0x_1\cdots x_{t-1})  (x_1x_2\cdots x_{t}) \cdots (x_{|\zeta(\omega_0)|-1}x_{|\zeta(\omega_0)|}\cdots x_{|\zeta(\omega_0)| +t-2}),
 $$
 which is the ordered list of the first $|\zeta(\omega_0)|$ words of length $t$ contained in $\zeta(\omega)$.
 
 Notice that the length of $\zeta_t(\omega)$ is equal to the length of $\zeta(\omega_0)$.  Thus, if $\zeta$ is of non-constant length, so is
 $\zeta_t$. 
 Let us look at the Fibonacci substitution which is of non constant length: $\zeta(0)=01, \zeta(1)=0$. There are three admissible words $00, 01, 10$ of length $2$, which constitute
 $\mathcal{A}_2$. 
  Since
  $$
    \zeta(00)=0101, \quad \zeta(01)=010, \quad \zeta(10)=001,
  $$
  we have 
  $$
     \zeta_2(00)=(01)(10), \quad  \zeta_2(01)=(01)(10), \quad  \zeta_2(10)=(00).
  $$
  Notice that $\zeta_2(10)$ has length $1$ because $\zeta(1)$ has length $1$. But $\zeta_2(00)$ and $\zeta_2(01)$ have length $2$.
  
  \medskip
 
 It is known that the primitivity of $\zeta$ implies the primitivity of $\zeta_t$ and the fix point of $\zeta_t$ is equal to
 $$
    U=(u_0u_1\cdots u_{t-1})(u_1u_2\cdots u_t) (u_2u_3\cdots u_{t+1})\cdots
 $$
 if $u=u_0 u_1 u_2\cdots$ is the fixed point of $\zeta$. So, the same idea as above permits the computation  of the frequencies of patterns of length $t$, namely the probability vector
 $\nu([\omega])_{\omega\in \mathcal{A}_t}$. It is the normalized Perron-Frobenius vector of the composition matrix $M_{\zeta_t}$ of $\zeta_t$. 
 \medskip

  Now we illustrate the method of computation by some concrete examples.  
  
  \subsection{Matrices selected by Fibonacci sequence}

Consider the Fibonacci sequence $\omega$, which is generated by the substitution $$
\zeta(0) =01, \quad \zeta(1) =0.
$$

\begin{cor} \label{cor:Fibonacci}Let $A_0 = \mathbf{u} \mathbf{v}'$ be a non-negative matrix of rank $1$ and $A_1=A=(a_{i,j})$ be a nonnegative matrix.
Let $\{\omega_n\}_{n \ge 0}$
be the Fibonacci sequence $01001010\cdots$. 
If $\mathbf{v}' \mathbf{u} >0$ 
and 
$\mathbf{v}'A\mathbf{u} >0$, we have
$$
 \lim_{n\to \infty} \frac{1}{n} \log \|A_{\omega_0}\cdots A_{\omega_{n-1}}\| =
 (\sqrt{5}-2) \log \mathbf{v}'\mathbf{u}  + \frac{3-\sqrt{5}}{2} \log \mathbf{v}'A \mathbf{u}. 
 $$
 If $\mathbf{v}' \mathbf{u} =0$ or
$\mathbf{v}'A\mathbf{u} =0$, then $A_{\omega_0}\cdots A_{\omega_{n-1}}=0$ when $n\ge 4$.
\end{cor}

\begin{proof} The last assertion is obvious because $A_0A_1A_0^2=0$ under the condition $\mathbf{v}' \mathbf{u} =0$ or
$\mathbf{v}'A\mathbf{u} =0$. In the sequel, we  compute the limit under the assumption $\mathbf{v}' \mathbf{u} >0$ and
$\mathbf{v}'A\mathbf{u} >0$.

We use the notation in Theorem \ref{thm:rank1}.  Since the word $11$ is forbidden, we have 
$\mathcal{R}=\{1\}$ and 
$F_1= \nu([1])=\nu([10])$.  We are led to compute the frequency $\nu([1])$.

Recall that the composition matrix of $\zeta$ is equal to 
$$
      M_{\zeta} = \begin{pmatrix}
       1 & 1 \\
       1&0\\
      \end{pmatrix}.
$$
Its Perron-Frobenius eigenvalue is $\frac{1}{1+\sqrt{5}}$ and the associated probability eigenvector is $(\frac{\sqrt{5}-1}{2}, \frac{3-\sqrt{5}}{2})$. 
Therefore  $\nu([0])=\frac{\sqrt{5}-1}{2}$ and $\nu([1])= \frac{3-\sqrt{5}}{2}$.
\end{proof}

\subsection{Matrices selected by Thue-Morse sequence}

 The Thue-Morse sequence $\omega$  is generated by the substitution $$
 \zeta: 0\mapsto 01, \quad 1\mapsto10.$$ 
 The Thue-Morse sequence is then
 $$0110100110010110\cdots.$$
 
  \begin{cor} Let $A_0 = \mathbf{u} \mathbf{v}'$ be a non-negative matrix of rank $1$ and $A_1=A=(a_{i,j})$ be a nonnegative matrix.
Let $\{\omega_n\}_{n \ge 0}$
be the Thue-Morse sequence. 
If $\mathbf{v}' \mathbf{u} >0$ 
and 
$\mathbf{v}'A\mathbf{u} >0$, we have
$$
 \lim \frac{1}{n} \log \|A_{\omega_0}\cdots A_{\omega_{n-1}}\| =
 \frac{1}{6} \Big(\log \mathbf{v}'\mathbf{u}  + \log \mathbf{v}'A \mathbf{u} + \log \mathbf{v}'A^2 \mathbf{u}  \Big).
 $$
 If $\mathbf{v}' \mathbf{u} =0$ or
$\mathbf{v}'A\mathbf{u} =0$, then $A_{\omega_0}\cdots A_{\omega_{n-1}}=0$ when $n\ge 4$.
\end{cor}

\begin{proof}   The Thue-Morse sequence does not contain the cube $1^3=111$, but contains $1$ and $11$. Thus  $\mathcal{R}=\{1, 11\}$.
 In order to apply Theorem \ref{thm:positive}, 
 let us first compute the frequencies of $2$-words. Observe that
 $$
   \zeta_2(00) = (01)(10), \ \    \zeta_2(01) = (01)(11), \ \   \zeta_2(10) = (10)(00), \ \   \zeta_2(11) = (10)(01).
 $$
 The composition matrix of $\zeta_2$ is then equal to 
 \begin{equation}\label{eq:M2-TM}
     M_{\zeta_2} = \begin{pmatrix}
         0 & 0 & 1 & 0\\
           1 & 1 & 0 & 1\\
           1 & 0 & 1 & 1\\
           0 & 1 & 0 & 0
     \end{pmatrix}.
 \end{equation}
 Its Perron-Frobenius eigenvalue is $2$ and the associated probability eigenvector is 
 $(\frac 16,\frac 13,\frac 13,\frac 16)$.  
We have $F_{11}=\nu([11])= \frac{1}{6}$. By the formula (\ref{eq:Fw}) in Theorem \ref{thm:ExactF}, we get
$$
F_{1} = \nu([1])- 2\nu([11])= [\nu([1])- \nu([11])] - \nu([11])= \nu([10])-\mu([11]) = \frac{1}{3}-\frac{1}{6}=\frac{1}{6}.
$$
It follows that 
$
    1 - (2F_1 +3F_{11}) =  \frac{1}{6}.
$
Thus the formula (\ref{eq:LiapRank2}) applies and gives what we want. 
\end{proof}

\subsection{Computation of exact frequencies through the normalized derivative substitution}

We have presented above the method, due to P. Michel \cite{Michel1974},  using the induced substitution to compute the exact frequencies. 
If there are return words of length $t$, the
$\zeta_t$-matrix of the induced substitution
$\zeta_t$ is a $|\mathcal{A}|^t \times |\mathcal{A}|^t$ matrix, where $\mathcal{A}$
is the alphabet on which the substitution is defined. For example, if we consider 
a substitutive sequence of two letters which has a return word of length $10$, then the matrix size is 
$2^{10}\times 2^{10}$.
When $t$ is larger, it is a large matrix and it is consuming for numerical calculation. In this way, we compute the frequencies of all words of length $t$.

However we do not need all frequencies but only some of them. 
There is another method, due to Fabien Durand \cite{Durand1998} (cf. also \cite[Chapter 11]{BDP2024}). As pointed out in \cite{Durand2011}, the first observation was made by G\'erard Rauzy,  
using another induced substitution called the derivative substitution.
If there are $R$ return words, the substitutive matrix of the derived substitution is of size $R\times R$. 
For example, if we consider 
a substitutive sequence of two letters which has  $100$ return words, then the size of involved matrix  is 
$100\times 100$. This is to be compared with the above mentioned matrix of size $2^{10}\times 2^{10}$, i.e. $1024 \times 1024$. 

Let us briefly recall the derived substitution method adapted to our setting. 
Consider a primitive substitution $\zeta$ on the alphabet 
$\mathcal{A}:=\{0,1, \cdots, s-1\}$ ($s\ge 2$). Let $\omega$ be a fix point of $\zeta$. Recall that $(X_\zeta, T)$ denotes the associated subshift. Now we consider the set $\mathcal{R}(0)$ of return words to the cylinder $[0]$: a word $\mathbf w$ belongs to $\mathcal{R}(0)$ if $\mathbf w0$ is a subword of $\omega$, and $\mathbf w$ contains $0$ as prefix and does not contain other $0$'s. In other words, $\mathbf w0=0*0$ where $*$ is a word without letter $0$ ($*$ can be an empty word).  Notice that $\mathcal{R}(0)$ is not the same as $\mathcal{R}$, which is defined earlier in this paper.
Also notice that it is possible that $0\in \mathcal{R}(0)$. 
Our substitutive sequence $\omega$ is clearly decomposed  in the the following way
\begin{equation}\label{eq:decom2}
    \omega = \mathbf w_00\mathbf w_10\mathbf w_2\cdots 0\mathbf w_n 0\cdots 
\end{equation}
where $\mathbf w_j \in \mathcal{R}(0)$ for all $j\ge 1$, and $\mathbf w_0$ is a word not containing the letter $0$, which can be the empty word.

Since $\zeta$ is primitive, the subshift $(X_\zeta, T)$ is minimal. Therefore the set of return words $\mathcal{R}(0)$ is finite, let $R= {\rm Card} \, \mathcal{R}(0)$. 
We relabel $\mathcal{R}(0)$ by $\mathcal{A}(0):=\{0,1, \cdots, R-1\}$. For  $i \in \mathcal{A}(0)$, we define $\phi(i)$ to be  the $i$-th return word, belonging to $\mathcal{R}(0)$, that appears in $\omega$. Thus $\phi$ defines a bijection between
$\mathcal{A}(0)$ and $\mathcal{R}(0)$, called the normalized return 
substitution. It defines a substitution $\eta$ on the alphabet
$\mathcal{A}(0)$ such that 
$$
     \phi \circ \eta = \zeta \circ \phi.
$$
This substitution $\eta$ is called the normalized derivative substitution of $\zeta$. Any bijection between $\mathcal{A}(0)$ and $\mathcal{R}(0)$ defines a substitution on $\mathcal{A}(0)$, called a derivative substitution.  However, we will only use the 
the normalized derivative substitution.

To illustrate these notions,  let us consider an example, the Thue-Morse substitution $\zeta$
defined by $0 \mapsto 01$, $1\mapsto 10$ and the 
corresponding Thue-Morse sequence $\omega =\zeta^\infty(0)$. 
We have $\mathcal{R}(0) =\{011, 01, 0\}$. The normalized return substitution is
$$
   \phi(0) =011, \quad \phi(1)= 01, \quad \phi(2) =0
$$
and the normalized derivative substitution:
$$
   \eta(0) =012, \quad \eta(1)=12, \quad \eta(2)= 1
$$
because
$$
    \zeta (011) = 011010=(011)(01)(0), \ \ \ 
    \zeta(01) =0110 =(011)(0), \ \ \ \zeta(0) =01 =(01).
$$
Therefore the $\eta$-matrix is given by
\begin{equation}\label{eq:M-TM}
 M_\eta= \begin{pmatrix}
      1 & 1& 0\\ 1 & 0& 1\\ 1 & 1 &0
  \end{pmatrix}.
\end{equation}

Let us return to the general case.  Let $\mu$ be the unique shift invariant measure on $X_\zeta$. Let $\tau$ be the first return time
to the cylinder $[0]$
$$
    \tau(x) =\inf \{n\ge 1: T^n x \in [0]\}.  
$$
The induced map $S: [0] \to [0]$,
where $[0]$ is understood as $[0]\cap X_\zeta$, is defined by $S(x) = T^{\tau(x)}(x)$. It is well known that 
$([0] , \mu_{[0]}, S)$ is ergodic, where $\mu_{[0]}$ is the normalized restriction of $\mu$ on 
$[0]$, namely
   $
      \mu_{[0]}(\cdot) = \frac{\mu(\cdot)\cap[0]}{\mu([0])}
   $.
   The  induced system $S$ is topologically conjugate to the substitutive system associated to $\eta$, which is primitive and thus is uniquely ergodic.  Therefore,
the ergodic system $([0], \mu_{[0]}, S)$  
is conjugate to the substitutive ergodic system
associated to the normalized derivative substitution $\eta$ (cf. \cite{DHS1999}, Lemma 19). 
Let $\lambda$  be the unique invariant measure of 
the system associated to 
the substitution $\eta$. We have $\mu_{[0]}=\lambda$.

Notice that the cylinder $[0]$ (in $X_\zeta$)
is decomposed into cylinders $[\mathbf u0]$ with $\mathbf u\in \mathcal{R}(0)$ so that 
    $$
    \mu([0]) = \sum_{\mathbf u\in \mathcal{R}(0)} \mu[\mathbf u0].
$$
For $\mathbf w\in \mathcal{R}$ (return word defined earlier), we have $0\mathbf w\in \mathcal{R}(0)$ and the exact frequency of $\mathbf w$ is equal to 
\begin{equation}\label{eq:EF2}
    F_{\mathbf w} =\mu([0\mathbf w0]) =  \mu_{[0]}([0\mathbf w0])\mu([0]) = \lambda([0\mathbf w]) \mu([0]).
\end{equation}

The formula \eqref{eq:EF2} allows us to compute the exact frequencies. Let us now look at some examples.
\medskip

{\bf Thue-Morse sequece.} For the Thue-Morse sequence 
 we have $\mathcal{R}=\{1,11\}$ but $\mathcal{R}(0)=\{011, 01, 0\}$. 
It is clear that the $\eta$-matrix (cf. \eqref{eq:M-TM}) has $1$ as Perron eigenvalue and $(\frac 13,
\frac 13,\frac 13)$ as the corresponding eigenvector. 
If follows
that $\lambda([011])=\lambda([01])=\lambda([0])=\frac{1}{3}$. As $\mu([0])=\frac{1}{2}$, by the formula
\eqref{eq:EF2}, we get $F_{11}=F_1= \frac{1}{6}$.  
This coincides with what we have computed using the
$\zeta_2$-matrix (cf. \eqref{eq:M2-TM}).
\medskip

{\bf Tribonacci sequence.}\ Now let us look at the Tribonnaci substitution:
   $$
    \zeta(0) =01, \quad \zeta (1) = 02, \qquad \zeta(2) = 0.
   $$
   Firstly we have
   \begin{equation}\label{eq:Trib}
 M_\zeta= \begin{pmatrix}
      1 & 1& 1\\ 1 & 0& 0\\ 0 & 1 &0
  \end{pmatrix}.
\end{equation}
The characteristic polynomial of $M_\zeta$ is $x^3-x^2-x-1$ and the Perron eigenvalue of $M_\zeta$
is equal to 
$$
\rho = \frac{1}{3} \left( \big( 19 + 3 \sqrt{33}\big) + \big( 19 - 3 \sqrt{33}\big) + 1\right) =1.839286755214161...
$$
The normalized Perron eigenvector is equal to 
$(\rho^{-1}, \rho^{-2}, \rho^{-3})$, so that 
$$
    \mu([0])= \frac{1}{\rho}, \quad  
    \mu([1])= \frac{1}{\rho^2}, \quad \mu([2])= \frac{1}{\rho^3}.
$$
Secondly, we have $\mathcal{R}=\{1, 2\}$, $\mathcal{R}(0)=\{01, 02, 0\}$ and 
$$\phi(0)=01,\  \phi(1)=02,\  \phi(2)=0.$$
It is observed that in this case, the normalized derivative substitution $\eta$ is equal to the initial Tribonacci substitution $\zeta$. This was first observed by Rauzy (cf. \cite{Durand2011}). By the formula \eqref{eq:EF2}, we get
$$
   F_{1} = \lambda([01])\mu([0])=\frac{1}{\rho}\cdot \frac{1}{\rho}= \frac{1}{\rho^2}, \qquad F_{2} = \lambda([02])\mu([0])= \frac{1}{\rho^2}\cdot \frac{1}{\rho}= \frac{1}{\rho^3}
$$
However, since $1$ and $2$ must be followed by $0$,
we can also argue that 
$$
F_1 =\mu([1])=\frac{1}{\rho^2}, \quad F_{2}=\mu([2])=\frac{1}{\rho^3}.
$$ 
Actually, for Tribonacci sequence, the exact frequencies of $1$ and $2$ are the frequencies of $1$ and $2$. Therefore, we do not actually need \eqref{eq:EF2}.

If we are given three non-negative matrix $A_0, A_1, A_2$ such that $A_0= \mathbf{u} \mathbf{v}'$
is of rank one and if we selection a sequence of 
matrices $\{A_{\omega_n}\}$ according to the
tribonacci sequence $\omega =(\omega_n)$.  Then,
by the formula \eqref{eq:LiapRankMain}, the Lyapunov exponent of $\{A_{\omega_n}\}$ is equal to 
$$
  L(\omega) = \left(\frac{1}{\rho}-\frac{1}{\rho^2} - \frac{1}{\rho^3} \right) \log \mathbf{v}' \mathbf{u} + \frac{1}{\rho^2} \log \log \mathbf{v}' A_1\mathbf{u}+   \frac{1}{\rho^3} \log \mathbf{v}' A_2\mathbf{u}.
$$

Similar results holds for the $m$-bonacci
sequence ($m\ge 4$) generated by the substitution
$$
0\mapsto 01,\ \  1\mapsto 02, \ \ \cdots, \ \ m-2\mapsto 0 (m-1),\ \  m-1 \mapsto 0. 
$$
\medskip
{\bf Another example.} Let us consider the substitution $\zeta(0)=01$, $\zeta(1) = 100110$. Then the substitutive sequence is decomposed into return words:
$$
(011)(0)(011)(01)(0)(011)(0)(01)(011)(0)(011)\cdots
$$
We have $\mathcal{R}(0)=\{011, 0, 01\}$ and 
$$
    011 \mapsto (011)(0)(011)(01)(0)(011)(0), \quad 0\mapsto (01), \quad 01\mapsto (011)(0)(011)(0).
$$
The $\zeta$-matrix and the $\eta$-matrix are respectively
$$
   M_\zeta =\begin{pmatrix}
       1 & 3\\ 1 &3
   \end{pmatrix}, \qquad 
   M_\eta =\begin{pmatrix}
       3 & 0& 2\\ 3&0&2\\ 1&1&0
   \end{pmatrix}.
$$
The matrix $M_\zeta$ has $2$ as Perron-eigenvalue with the eigenvector $(1/2, 1/2)$, and for the matrix $M_\eta$
the Perron eigenvalue is $4$ and its Perron
eigenvector (2/5, 2/5, 1/5). Therefore, by the formula
\eqref{eq:EF2}, we get the exact frequencies of $11$ and $1$:
$$
   F_{11} =\frac{2}{5}\cdot \frac{1}{2}=\frac{1}{5},
   \qquad 
   F_{1} =\frac{1}{5}\cdot \frac{1}{2}=\frac{1}{10}.
$$

If we use the method of Michel \cite{Michel1974}, we have to compute the  frequencies of words of length $3$. The matrix $M_{\zeta_3}$ is a $8\times 8$ matrix. 
This allows us to compute all the 
  frequencies of words of length $3$. But we only
  need that of $011$ and $010$. 
  
  With the help of the formula \eqref{eq:Fw}, we can only compute the  frequencies of words of length $2$.

\section{Matrices selected by $\mathcal{B}$-free integers}\label{sect:BF}

 In this section, we apply Theorem \ref{thm:rank1} and Theorem \ref{thm:ExactF} to compute the Lyapunov exponents of the products of matrices selected by the
 characteristic function of  $\mathcal{B}$-free integers. An interesting particular example is the  
 characteristic function of   square-free integers. An important fact is that the  characteristic function of   $\mathcal{B}$-free integers as a point in $\{0,1\}^\mathbb{N}$ is generic for the Mirsky measure and another important fact is that the distribution of Mirsky measure
 is well known (cf \cite{ALR2015}).

 \subsection{Characteristic function of $\mathcal{B}$-free integers and Mirsky measure}
 Following \cite{ALR2015}, we recall below some basic facts about the $\mathcal{B}$-free integers and the associated dynamical system, especially the Mirsky measure. 
 Let $$
      \mathcal{B}=\{b_k: k \ge 1\} \subset \{2,3, 4,\cdots\}
    $$ 
 be an increasing sequence of integers satisfying \\
 \indent (B1)   $b_k$ and $b_{k'}$ are relatively prime for any  $1\le k<k'$,\\
 \indent (B2)   $\sum_{k=1}^\infty \frac{1}{b_k}<\infty$.\\
 Natural numbers which have no factors in $\mathcal{B}$ are said to be $\mathcal{B}$-free.  
  The most interesting example is the set of squares of primes numbers and the corresponding $\mathcal{B}$-free integers are the square-free integers.
 The above notion of $\mathcal{B}$-free numbers was introduced by P. Erd\"{o}s \cite{Erdos1966} in 1966, who proved
 that there exists a constant $0<c<1$ such that each interval $[x, x+x^c]$ contain at least one $\mathcal{B}$-free number when $x$ is large enough.  There are many works trying to improve the constant $c$
 of Erd\"{o}s (cf. \cites{Wu1993,Wu1993b,KRW2007}). The fact that $c<1$ will be useful for us to apply Theorem \ref{thm:rank1}. See \cites{Mirsky1949,CS2013,Peckner2015} for studies on the classical square free integers. 
 
 The characteristic function of the set of $\mathcal{B}$-free integers is the sequence
 $\eta:=\eta_{\mathcal{B}}:= (\eta_n)_{n\ge 1} \in\{0,1\}^\mathbb{N}$ defined by 
 $$
     \eta_n = \left\{ \begin{array}{cc}   0   \  \  \  {\rm if } \ \ n \  {\rm is \  divided \ by \ some }\ b_k,\\
                                                  1 \  \  \  {\rm otherwise}. \hspace{7em} \
                            \end{array}    \right.
 $$ 
 Recall that the shift map $\sigma$ on $\{0,1\}^\mathbb{N}$ is defined,   for $x=(x_n) \in \{0,1\}^\mathbb{N}$,  by $\sigma x =(x_n')$ where $x'_n= x_{n+1}$ for all $n\ge 0$. The closure $X_\eta:=\overline{\{\sigma^n\eta: n\ge 0\} }$ of the orbit of $\eta$ under the shift is a closed set such that $\sigma(X_\eta)\subset X_\eta$. 
 We are interested in the dynamical properties of the subshift $(X_\eta, \sigma)$. 
 
 The subshift $X_\eta$ is well described by $\mathcal{B}$-admissible sequences (cf. \cite{ALR2015}). A subset $A \subset \mathbb{N}$ is said to be   {\bf $\mathcal{B}$-admissible} if 
 $t(A, b_k)<b_k$ for all $k\ge 1$ where  $t(A, b)$ is the number of classes modulo $b$ in $A$, i.e.   
 $$
   t(A, b)  = \#\{  a \!\!\mod b : a \in A\}.
 $$
 The {\bf support} of a sequence $x =(x_n)_{n\ge 1} \in \{0,1\}^\mathbb{N}$ is the set $\{n\in \mathbb{N}: x_n=1\}$, denoted ${\rm supp}\, x$. 
 An infinite  sequence $x\in \{0,1\}^\mathbb{N}$ is said to be  $\mathcal{B}$-admissible if its support is $\mathcal{B}$-admissible.  In the same way, a finite sequence
 $x_1x_2\cdots x_N \in \{0,1\}^N$
 is said to be  $\mathcal{B}$-admissible  if its support $\{n \in \{1,2,\cdots, N\}: x_n=1\}$ is $\mathcal{B}$-admissible. Let $X_{\mathcal{B}}$ be the set of all  $\mathcal{B}$-admissible sequences.
 It was proved that $X_\eta= X_{\mathcal{B}}$, and 
 it was also proved that  the sequence $\eta$ is  $\nu_{\mathcal{B}}$-generic for some invariant measure $\nu_\mathcal{B}$, called Mirsky measure (\cite{ALR2015}).
 
 Let us get together the above mentioned properties and other properties 
 of the subshift $(X_\eta, \sigma)$ in the following two  theorems.
 
 \begin{thm} [\cite{ALR2015}]\label{thm:ALR1}
 The following are true:\\
 \indent {\rm (1)}  $X_\eta= X_{\mathcal{B}}$.\\
 \indent {\rm (2)} The topological entropy of  $(X_\eta, \sigma)$ is equal to $\prod_{k=1}^\infty (1- 1/b_k)$.\\
 \indent {\rm (3)}  The subshift $X_\eta$  supports a $\sigma$-invariant measure $\nu_{\mathcal{B}}$ and the measure preserving system $(X_\eta, \nu_{\mathcal{B}}, \sigma)$ is isomorphic to the minimal and unique ergodic rotation by adding $(1,1,\cdots)$ on the compact abelian group   
 $\prod_{k=1}^\infty \mathbb{Z}/b_k\mathbb{Z}$.\\
  \indent {\rm (4)}  $\eta$ is $\nu_{\mathcal{B}}$-generic.
 \end{thm} 
 
 Let $T:  \prod_{k=1}^\infty\mathbb{Z}/b_k\mathbb{Z} \to  \prod_{k=1}^\infty \mathbb{Z}/b_k\mathbb{Z}$  be the rotation defined by
 $T \omega = (\omega_k +1)_{k\ge 1}$. Notice that here we add $1$ to every coordinate and, thus, $T$ is  not the odometer transformation.
 
  Let us define the {\bf $\mathcal{B}$-free test function} on $\prod_{k=1}^\infty\mathbb{Z}/b_k\mathbb{Z}$:
 $$
   B(\omega) =  \left\{ \begin{array}{cc}   0   \  \  \  {\rm if } \ \ \omega_k =0 \! \mod b_k \ \  {\rm for \ some }\ k\ge 1,\\
                                                  1 \  \  \   {\rm if } \ \ \omega_k \not =0 \!  \mod b_k \ \  {\rm for \ all }\ k\ge 1. \ \ \ 
                            \end{array}    \right.
 $$
 For $n\ge 1$, define $$\underline{n} = (n\!\! \mod b_1,  n\!\! \mod b_2, \cdots)\in \prod_{k=1}^\infty\mathbb{Z}/b_k\mathbb{Z}.$$ 
 It is clear that $n$ is $\mathcal{B}$-free if and only if $B(\underline{n})=1$. So, every integer $n$ is identified with   a point $\underline{n}$
 in $\prod_{k=1}^\infty\mathbb{Z}/b_k\mathbb{Z}$ (the map $n \mapsto \underline{n}$ being injective but not surjective) and  the value $B(\underline{n})$ indicates whether $n$ is $\mathcal{B}$-free or not. 
 As $$
 \underline{n}= (n \!\!\mod b_1, n \!\!\mod b_2, \cdots)= T^n(\underline{0}),
 $$  we get 
 $$\eta=(B(T^1(\underline{0})), B(T^2(\underline{0})), B(T^3(\underline{0})), \cdots ). 
 $$
 Define the map $\varphi: \prod_{k=1}^\infty\mathbb{Z}/b_k\mathbb{Z} \to X_{\mathcal{B}}\subset \{0,1\}^\mathbb{N}$ by 
 $$
     \varphi(\omega) = (B(T^n \omega))_{n\ge 1}, \ \ \ {\rm i.e.}\ \ \   \varphi(\omega)  =(B(\omega +\underline{n}))_{n\ge 1}.
 $$
 Indeed, $\varphi(\omega)$ is $\mathcal{B}$-admissible because for any $k$,
  $$
  {\rm supp}\, \varphi(\omega)\!\! \mod b_k
  =\{n \!\!\mod b_k: B(\omega +\underline{n})=1\},$$ 
  which does not contain $-\omega_k$.

 The map $\varphi$ is not continuous but measurable and $\sigma \circ T=T\circ \varphi$.
 Let $\mathbb{P}$ be the normalized Haar measure on  $\prod_{k=1}^\infty\mathbb{Z}/b_k\mathbb{Z}$. The {\bf Mirsky measure}
 $\nu_{\mathcal{B}}$
 is defined to be the image measure of $\mathbb{P}$ under $\varphi$. 
 
 For two disjoint finite subset $A, B \subset \mathbb{N}$ we define the cylinder set
 $$
   C_{A, B} = \{(x_n)_{n\ge 1}\in \{0,1\}^\mathbb{N}: \forall n \in A, x_n=1; \forall m \in B, x_m=0\}.
 $$
 We write $C^1_A$ for $C_{A, \emptyset}$, $C^0_B$ for $C_{\emptyset, B}$. If $A=\{1,2, \cdots, t\}$, we usually 
 write $[1^t]$ for $C_A^1$. 
 
  \begin{thm} [\cite{ALR2015}]\label{thm:ALR2}
  The following are true.\\
 \indent {\rm (1)} $(X_\eta, \nu_{\mathcal{B}}, \sigma)$ is a factor of  $(\prod_{k=1}^\infty\mathbb{Z}/b_k\mathbb{Z}, \mathbb{P}, T)$. So, $\nu_{\mathcal{B}}$ is ergodic and of zero entropy.\\
 \indent {\rm (2)} For finite subset $A\subset \mathbb{N}$ we have
 $$
        \nu_{\mathcal{B}}(C^1_A) = \prod_{k=1}^\infty \left(1- \frac{t(A, b_k)}{b_k}\right).
 $$ 
 \indent {\rm (3)} For finite disjoint subsets $A,B\subset \mathbb{N}$ we have
 $$
        \nu_{\mathcal{B}}(C_{A, B}) = \sum_{A\subset D\subset A\cup B} \prod_{k=1}^\infty (-1)^{\#(D\setminus A)}\left(1- \frac{t(A, b_k)}{b_k}\right).
 $$ 
 \indent {\rm (4)} $A$ is $\mathcal{B}$-admissible iff $\nu_{\mathcal{B}}(C^1_A)>0$ iff $\nu_{\mathcal{B}}(C_{A,B})>0$ for all $B$ disjoint from $A$.
 \end{thm}

 \subsection{Lyapunov exponent of the product of the matrices $\{A_{\eta_n}\}$}
 Given two  non-negative  matrices $A_0$ and $A_1$  with $A_0$ having rank $1$, we are going to compute the Lyapunov exponent of the product of the matrices
 $\{A_{\eta_n}\}$ selected by the characteristic function $\eta=(\eta_n)$ of the sequence of $\mathcal{B}$-free integers.
 
  Recall that  $
      \mathcal{B}=\{b_k: k \ge 1\} \subset \{2,3, 4,\cdots\}
    $ 
 be an increasing sequence of integers such that 
  $b_k$ and $b_{k'}$ ($k\not=k'$) are relatively prime and
 $\sum_{k=1}^\infty \frac{1}{b_k}<\infty$. 
 
 For $1\le t<b_1$, we introduce the function of Euler type
 $$
     \zeta^{\mathcal{B}}_t(s) = \prod_{k=1}^\infty \left( 1-\frac{t}{b_k^{s/2}}\right)^{-1}.
 $$
 When $\mathcal{B}=\{p^2: p\in \mathcal{P}\}$ ($\mathcal{P}$ being the set of prime numbers),  $ \zeta^{\mathcal{B}}_1(s) $ is the Riemann zeta function
 $\zeta(s) =\sum_{n=1}^\infty \frac{1}{n^s}$.   Clearly we have  $1<\zeta^{\mathcal{B}}_t(2)<\infty$ for $1\le t<b_1$.
 
 \begin{thm} \label{thm:MA_Bfree}
 Let  $\eta$ be the characteristic function of the set of $\mathcal{B}$-free integers.
 Let $A_0$ and $A_1$ be two non-negative matrices. Suppose $A_0 = \mathbf{u} \mathbf{v}'$ is of rang one
 and that $\mathbf{v}'\mathbf{u}\not=0$, $\mathbf{v}'A_{\mathbf{w}}\mathbf{u}\not=0$ for $\mathbf{w} \in 
 \mathcal{R}=\{1, 11, \cdots, 1^{b_1-1}\}$.
 Then 
 the Lyapunov exponent  of the sequence of matrices $(A_{\eta_n})$ is equal to 
 \begin{equation}\label{eq:LiapBfree}
    L = \left(1- \sum_{\mathbf{w}\in \mathcal{R}} (|\mathbf{w}|+1) F_{\mathbf{w}} \right) \log \mathbf{v}'\mathbf{u} + \sum_{\mathbf{w}\in \mathcal{R}} F_{\mathbf{w}}\log \mathbf{v}'A_{\mathbf{w}}\mathbf{u},
 \end{equation}
 where 
 \begin{equation}\label{eq:FwB}
      F_{\mathbf{w}} = \sum_{j=0}^{b_1-1 -|\mathbf{w}|} (-1)^j S^{(j)}_{\mathbf{w}} 
        \end{equation}
        for $\mathbf{w}\in \mathcal{R}$, with $S^{(0)}_{\mathbf{w}} =\frac{1}{\zeta^{\mathcal{B}}_{|\mathbf{w}|}(2)}$ and for $1\le s<b_1-1-|\mathbf{w}|$
  \begin{equation}\label{eq:3*}
       S_{\mathbf{w}}^{(s)} =
    \sum_{j_s=s}^{b_1-1-|\mathbf{w}|} \frac{1}{\zeta^{\mathcal{B}}_{|\mathbf{w}|+j_s}(2)}  \sum_{1\le j_1<j_2<\cdots < j_{s-1}<j_s}(j_1+1) (j_2-j_1+1) \cdots (j_s-j_{s-1} +1).
 \end{equation}
 \end{thm} 
 
 \begin{proof} The subset  $\{1,2, \cdots, b_1\}$ is not $\mathcal{B}$-admissible as it is a complete system modulo $b_1$. Therefore there is no pattern $1^t$ in $\eta$ for $t\ge b_1$. 
 But all proper subsets of  $\{1,2, \cdots, b_1\}$ are clearly $\mathcal{B}$-admissible. Therefore we get the set of return words of $\eta$ (return to $0$)
 $$
    \mathcal{R}=\{1, 11, \cdots, 1^{b_1-1}\}.
 $$
 
 Theorem \ref{thm:rank1} is applicable:    firstly,  because $\mathcal{R}$ is finite and the condition 
 (iv) in Theorem \ref{thm:rank1} is also satisfied for the reason of  existence of $\mathcal{B}$-free integers in small interval and the  result  of Erd\"{o}s mentioned in the introduction is sufficient  (cf. \cites{Wu1993, Wu1993b,KRW2007}
 for better results on the existence of $\mathcal{B}$-free integers in small intervals).
 Secondly,  because $\eta$ is $\nu_\mathcal{B}$-generic, by Theorem \ref{thm:ALR1} (4). The formulas (\ref{eq:LiapBfree}) and (\ref{eq:FwB}) follows directly from Theorem \ref{thm:rank1}.
 
 What we need is just the computation of $S_{\mathbf{w}}^{(s)}$ for  $\mathbf{w}=1^k \in \mathcal{R}$ ($1\le k <b_1-1$), which is relatively straightforward. First,  
 $$
     S_{1^k}^{(1)} =\sum_{j=1}^{b_1-1-k} N_{1^k}(1^{k+j}) \nu_{\mathcal{B}}([1^{k+j}]) = \sum_{j=1}^{b_1-1-k} (j+1) \nu_{\mathcal{B}}([1^{k+j}]),
 $$
 just because $N_{1^k}(1^{k+j}) = j+1$. For the same reason, we have
 \begin{eqnarray*}
     S_{1^k}^{(2)} & = & \sum_{1\le j_1<j_2\le b_1-1-k} N_{1^k}(1^{k+j_1}) N_{1^{k+j_1}}(1^{k+j_2}) \nu_{\mathcal{B}}([1^{k+j_2}]) \\
     & = & \sum_{1\le j_1<j_2\le b_1-1-k} (j_1+1) (j_2-j_1+1) \nu_{\mathcal{B}}([1^{k+j_2}])\\
        & = & \sum_{j_2=2}^{b_1-1-k} \nu_{\mathcal{B}}([1^{k+j_2}])  \sum_{j_1=1}^{j_2-1}(j_1+1) (j_2-j_1+1) .
 \end{eqnarray*}
 In general, for $1\le s < b_1-1-k$, we have 
  \begin{eqnarray*}
     S_{1^k}^{(s)} & = &  \sum_{1\le j_1<j_2<\cdots <j_s \le b_1-1-k} N_{1^k}(1^{k+j_1}) N_{1^{k+j_1}}(1^{k+j_2}) \cdots N_{k+j_{s-1}}(1^{k+j_s}) \nu_{\mathcal{B}}([1^{k+j_s}]) \\
  & = & \sum_{j_s=s}^{b_1-1-k} \nu_{\mathcal{B}}([1^{k+j_s}])  \sum_{1\le j_1<j_2<\cdots j_{s-1}<j_s}(j_1+1) (j_2-j_1+1) \cdots (j_s-j_{s-1} +1).
  \end{eqnarray*}
This is (\ref{eq:3*}) because,  by Theorem \ref{thm:ALR2} (2),  we have
 $$
         \nu_\mathcal{B}([1^t]) = \frac{1}{\zeta^\mathcal{B}_t(2)} \ \ \ {\rm for \ all}  \ \  \  1\le t\le b_1-1. 
 $$

 \end{proof}

 Remark that $\mathbf{w} \in \mathcal{R}$ takes the form $1^t$ with $1\le t<b_1$, for which  $A_{\mathbf{w}}= A_1^t$. 
 
 \subsection{Square-free integers}
 
 Let us look at a concrete example.
 In the case of square-free integers, we have $b_1=4$ and $\mathcal{R}=\{1,11, 111\}$. For simplicity, we write $\nu$ for $\nu_{\mathcal{B}}$. We have seen in { Example 2}, Subsection 2.2, that
 $$
    F_{111} = \nu([111]), \quad F_{11}= \nu([11]) -2 \nu([1]), \quad F_1 = \nu([1]) - 2 \nu([11])  + \nu([111]).
 $$
 
 We are going to give  another way to compute the exact frequencies $F_{1}, F_{11}, F_{111}$. 
 Let us first examine the patterns of the forms $1^t$ and $0^s$ in $\eta$, or the patterns of the forms $01^t0$ and $10^s1$, which will help us to compute both the exact frequencies of  exactly $t$ consecutive $1$
 or  $s$ consecutive $0$'s (the exact frequency  of $s$ consecutive $0$'s will be not used).  
 
 First notice that patterns $01^t0$ with $t\ge 4$ are impossible, because among four consecutive integers $n, n+1,n+2,n+3$ there is one which is congruent to $0$ modulus $4$.
However patterns $10^s1$ are possible for all $s\ge 1$. 
 
 Let $u_1\cdots u_\ell \in \{0, 1\}^\ell$,  the frequency of $u_1\cdots u_\ell$ is given  by  
 $$
  \nu([u_1\cdots u_\ell]) =\lim_{N\to \infty} \frac{1}{N}\sum_{n=0}^{N-1} \mathbf{1}_{[u_1\cdots u_\ell]}(\sigma^n \eta).
 $$
 Notice that the exact frequencies $F_{1^t}$ of $1^t$ ($1\le t\le 3$) and $F_{0^s}$ of  $0^s$ ($s\ge 1$) are respectively  the exact frequencies of $1^t0$ and of  $0^s1$. In other words,
$F_{1^t}=F_{1^t 0}$ and $F_{0^s}=F_{0^s1}$. 
 
 In each pattern $1^{t'}0$ with $t'\ge t$ we see exactly one occurrence of $1^t 0$. Similarly,  in each pattern $0^{s'}1$ with $s'\ge 1$ we see exactly one occurrence of $0^s 1$.
 Thus we have the following relations
 $$
      \nu([1^t 0]) = F_{1^t 0} +  F_{1^{t+1} 0} + F_{1^{t+2} 0} + \cdots; \qquad   \nu([0^s 1]) = F_{0^s 1} +  F_{0^{s+1} 1} + F_{0^{s+2} 1} + \cdots.
 $$
 It follows that 
 \begin{equation} \label{Ft}
        F_{1^t}= F_{1^t 0} =\nu([1^t 0]) - \nu([1^{t+1} 0]) ; 
 \end{equation}
  \begin{equation}\label{Fs}
         F_{0^s} =F_{0^s 1} = \nu([0^s 1]) - \nu([0^{s+1} 1]).
 \end{equation}
 We are led to compute $\nu([1^t 0])$ and  $\nu([0^s 1])$.
 
  For $a \in \mathbb{R}$, we introduce the function 
 $$
 \zeta_a(s) = \prod \Big(1- \frac{a}{p^{s}}\Big)^{-1} =   \sum_{n=1}^\infty \frac{a^{\Omega(n)}}{n^{s}} , 
 $$
 where the product is taken over all prime numbers and $\Omega(n)$ is the number of primes (counted with multiplicity) contained in $n$.
 The function $\zeta_a(s)$ is analytic in ${\rm Re}\, s>1$ and $\zeta_1(s)$ is equal to the Riemann zeta function $\zeta(s)$. 
 We have $$
 \zeta_3(2) = 7.968954\cdots >\zeta_2(2)=3.099486\cdots  > \zeta_1(2) =\frac{\pi^2}{6}=1.644934\cdots.
 $$
 
 Recall that for $t\ge 4$, $F_{1^t 0} =\nu([1^t 0])- \nu([1^{t+1} 0])=0$.  By Theorem \ref{thm:ALR2} (3), we have
 \begin{eqnarray*}
     \nu([10]) & =& \prod \left(1- \frac{1}{p^2}\right) - \prod \left(1- \frac{2}{p^2}\right) =\frac{1}{\zeta_1(2)} -\frac{1}{\zeta_2(2)} = 0.285293,\\
     \nu([110]) & =& \prod \left(1- \frac{2}{p^2}\right) - \prod \left(1- \frac{3}{p^2}\right) = \frac{1}{\zeta_2(2)} -\frac{1}{\zeta_3(2)}= 0.197147,\\
     \nu([1110]) & =& \prod \left(1- \frac{3}{p^2}\right) - \prod \left(1- \frac{4}{p^2}\right) = \frac{1}{\zeta_3(2)}=0.125487.
 \end{eqnarray*}
 (the last product is equal to $0$). So, by (\ref{Ft}) we get
 \begin{eqnarray*}
     F_1=F_{10} & =&\frac{1}{\zeta_1(2)} -\frac{2}{\zeta_2(2)} + \frac{1}{\zeta_3(2)}=0.0881459;\\
     F_{11}=F_{110} & =&\frac{1}{\zeta_2(2)} - \frac{2}{\zeta_3(2)}=0.0716601;\\
    F_{111}= F_{1110} & =& \frac{1}{\zeta_3(2)}=0.125487.
 \end{eqnarray*}

 We can similarly find $F_{0^s1}$, which will be not necessary for us to compute the Lyapunov exponent. In order to find all $F_{0^s 1}$, we are led to compute $\nu([0^r 1])$ for $r\ge 1$, by (\ref{Fs}).  Notice that 
 $$
    \nu([0^r 1]) = \nu([0^r]) - \nu([0^{r+1}]).
 $$
 By Theorem \ref{thm:ALR2} (2), we can compute $ \nu([0^r])$ by 
 \begin{equation}
     \nu([0^r]) = \sum_{D\subset \{1,2, \cdots, r\}} (-1)^{r-|D|} \prod \left( 1- \frac{t(D, p^2)}{p^2}\right)
 \end{equation}
 where $t(D, p^2) =\#\{d \mod p^2: d \in D\}$.
 \medskip

 Let us look at the patterns $1, 11, 111$ among the first 300 terms of $\eta$ (50 by row, grouped by 10, the first of each group is marked blue):
 $$
  {\color{blue}1} 1 1 0111001 {\color{blue}1} 00110101 0  {\color{blue}1} 11001001 1{\color{blue}1} 01110111 0{\color{blue}1} 110001000
  $$
  $$
   {\color{blue}1} 010101110 {\color{blue}1} 100111011 {\color{blue}1} 011001110 {\color{blue}0} 110111010 {\color{blue}1} 011101000
 $$
 $$
     {\color{blue}1} 110111011  {\color{blue}1} 0111 00110  {\color{blue}0} 110001011 {\color{blue}1} 011101110  {\color{blue}1} 1101 10010   
 $$
$$
    {\color{blue}1} 001101110 {\color{blue}1} 010111001 {\color{blue}0}011001110 {\color{blue}1}  110111001  {\color{blue}1} 011101010
$$
$$
   {\color{blue}1} 110110011 {\color{blue}1} 011101110 {\color{blue}1} 110011011 {\color{blue}1} 010101110 {\color{blue}1} 000011101
   $$
   $${\color{blue}1} 011101110 {\color{blue}0} 110111010 {\color{blue}1} 011001100 {\color{blue}1} 110111001 {\color{blue}1} 0101 00010
   $$
   The pattern $0^4$ is seen only once for $\eta(242) \eta(243)\eta(244) \eta(245)$. Recall that
   $$
    242=2 \times 11^2 , \quad 243= 3^5, \quad 244=2^2\times 61, \quad 245=5\times 7^2
   $$
   all of which contain a square factor.
   \medskip

  \section{Multifractal analysis}\label{sect:MA}

 Here we present an application of Theorem  \ref{thm:rank1} and Theorem \ref{thm:ExactF}.  to the computation of multifractal spectrum of weighted Birkhoff averages. 
 
 Let us first recall a theoretical result from  \cite{Fan2021}, see also \cites{BRS2021,BRS2022} for information on the same topic.
 Let $S$ be a finite set of cardinality $|S|$. 
We consider the space $X=S^\mathbb{N}$ on which is defined the left-shift map $T$, a continuous function $f: X\to \mathbb{R}$ and a weight sequence $\{w_k\} $ taking values in $\mathbb{R}$.  We are interested in the multifractal analysis  of the wighted ergodic average 
$$
   Mf(x) := \lim_{n\to \infty}\frac{1}{n}\sum_{k=0}^{n-1} w_k f(T^k x).
$$ 
Following \cite{Fan2021},  we first define the pressure function
$$
   \psi(\beta) = \lim_{n\to \infty} \frac{1}{n} \log \int_X \exp\left( \beta \sum_{k=0}^{n-1} w_k f(T^k x)\right) dx +\log |S|,
$$
if the limit exists, where $\beta \in \mathbb{R}$ and $dx$ denotes the symmetric Bernoulli measure on $S^\mathbb{N}$. When $w_k=1$ for all $k$, we recover the classical pressure function.
It is clear that $\psi$ is a convex function.
For $\alpha \in \mathbb{R}$, we define the level set
$$
     E(\alpha) = \{x\in X: Mf(x) =\alpha\}.
$$
Recall that the usual metric on $X$ is defined by $d(x,y) = |S|^{-n}$ where $n = \inf \{k\ge 0: x_k\not=y_k\}$. The Hausdorff dimension $\dim$  and the packing dimension  
$\text{Dim}$ in the following theorem refer to this metric.

\begin{thm}[\cite{Fan2021}, Theorem 1.1] Assume that $f$ is a H\"{o}lder continuous function. 
If the pressure function $\psi$ is well defined  and differentiable on $\mathbb{R}$, then for $\alpha = \psi'(\beta)$ we have
\begin{equation}\label{eq:Mspect}
    \dim E(\alpha) =\Dim\, E(\alpha) = \frac{\psi(\beta)- \alpha \beta}{\log |S|} =  - \frac{\psi^*(\alpha)}{\log |S|} 
\end{equation}
where 
$\psi^*$ denotes the Legendre transform of $\psi$.
\end{thm}

Thus the multifractal analysis is reduced to the computation of the pressure function $\psi$ and its Legendre transform
$$
   \psi^*(\alpha) = \sup_\gamma (\alpha \gamma-\psi(\gamma)).
$$ 

Now assume that $f(x)=f(x_0, x_1)$ depends only on the first two coordinates (when $f$ depends on a finite number of coordinates, we can reduce it to the case of dependence on the first two coordinates). We introduce the following positive matrices 
$$
 A_{w_k}(\beta):=(e^{\beta w_k f(i, j)})_{(i, j)\in S\times S}.
$$
 It is easy to see (cf. \cite{Fan2021}, p.30) that  the limit defining $\psi(\beta)$ exits if and only if the following Lyapunov exponent $L(\beta)$ exists:
$$
   L(\beta) = \lim_{n\to \infty} \frac{1}{n} \log \|A_{w_0}(\beta) A_{w_1}(\beta)   \cdots A_{w_{n-1}}(\beta)  \|.
$$
In this case,  we have
 $\psi(\beta)=L(\beta)$.
Therefore we can compute $\psi$ by applying Theorem  \ref{thm:positive}.  
In the following, we are going to illustrate  how it works by examples. 

\subsection{Fibonacci weight}
    Consider the Fibonacci sequence generated by the substitution $0\mapsto 01$ and $1\mapsto 0$, which is a typical primitive substitutive sequence. Consider a function $f(x_0,x_1)$ on $\{0, 1\}^2$.
For $\beta \in \mathbb{R}$, let 
\begin{eqnarray*}
    A_0(\beta) & = & \begin{pmatrix}
      e^{0\cdot f(0,0)} & e^{0\cdot  f(0,1)} \\ e^{ 0\cdot f(1,0)} & e^{0\cdot f(1,1)}  
     \end{pmatrix}
    = \mathbf{1} \mathbf{1}', \\
     A_1(\beta)& = &  \begin{pmatrix}
      e^{\beta \cdot f(0,1)} & e^{\beta\cdot  f(0,1)} \\ e^{ \beta \cdot f(1,0)} & e^{\beta \cdot f(1,1)}  
     \end{pmatrix}.
\end{eqnarray*}
In this case, we have 
$$\mathbf{1}'\mathbf{1}=2, \qquad \mathbf{1}'A \mathbf{1}= e^{\beta f(0,0)} +e^{\beta f(0,1)} + e^{\beta f(1,0)} +e^{\beta f(1,1)}.
$$ So, 
by Corollary \ref{cor:Fibonacci}, the Lyapunov exponent 
$$
\psi(\beta):=\lim_{n\to \infty} \frac{1}{n} \log \|A_{\omega_1}(\beta)\cdots A_{\omega_n}(\beta)\| 
$$
is well defined and is equal to 
\begin{equation}\label{eq:Fibonacci_2x2}
\psi(\beta)=  (\sqrt{5}-2)\log 2 + \frac{3-\sqrt{5}}{2} \log (e^{\beta f(0,0)} +e^{\beta f(0,1)} + e^{\beta f(1,0)} +e^{\beta f(1,1)}).
\end{equation}
Clearly we have
$$
    \psi'(\beta) = \frac{3-\sqrt{5}}{2} \cdot \frac{f(0,0)e^{\beta f(0,0)} + f(0,1)e^{\beta f(0,1)} + f(1,0)e^{\beta f(1,0)} + f(1,1) e^{\beta f(1,1)}}{e^{\beta f(0,0)} +e^{\beta f(0,1)} + e^{\beta f(1,0)} +e^{\beta f(1,1)}}
$$
and 
$$
   \psi'(+\infty) = \frac{3-\sqrt{5}}{2} \cdot \max_{i, j} f(i, j), \qquad   \psi'(-\infty) = \frac{3-\sqrt{5}}{2} \cdot \min_{i, j} f(i, j).
$$
The multifractal spectrum is supported by the interval $[\psi'(-\infty), \psi'(+\infty)]$. See Figure \ref{fig:Fibonacci} for the graphs of $\psi$
and of the multifractal spectrum in the special case $f(x_0, x_1) = x_0x_1$. 
 \begin{figure}[htb]
 	\centering
 	\includegraphics[width=0.45\linewidth]{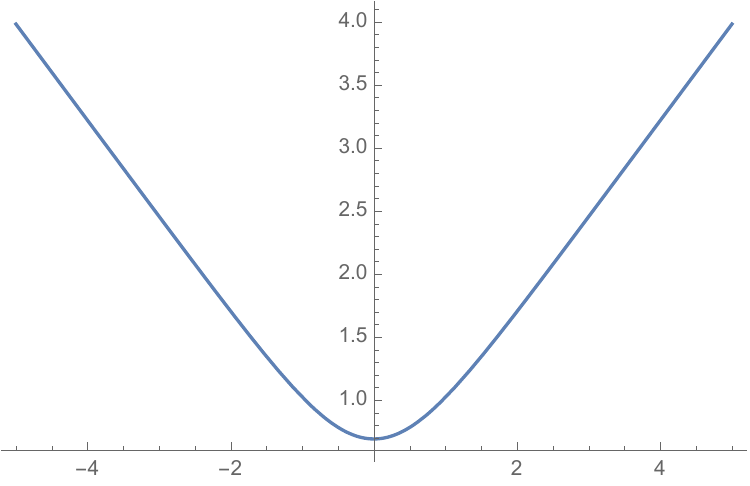}
 	\includegraphics[width=0.45\linewidth]{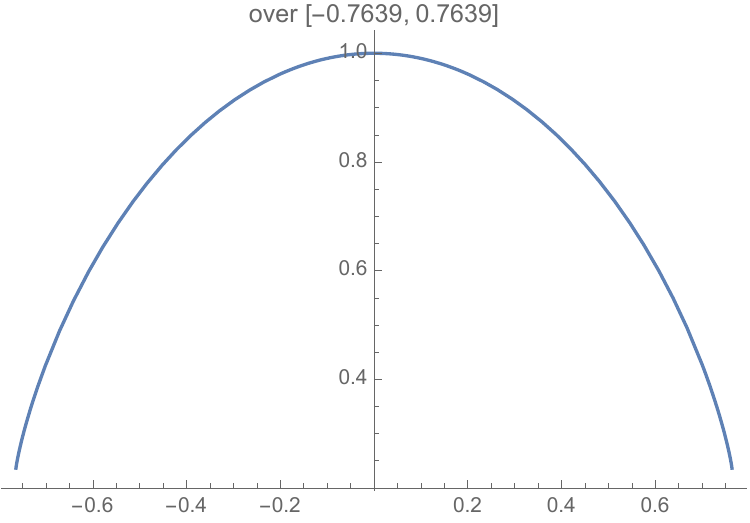}
 	\caption{ Graphs of  $\psi(\beta)$ and multifractal spectrum in the case of   $f(x_0, x_1) = x_0x_1$ with Fibonacci weight.}
 	\label{fig:Fibonacci}
 \end{figure}

\medskip

 \subsection{M\"obius weight}
 We can  determine the multifractal spectrum of the Birkhoff average
  $$
      Mf(x):=  \lim_{N\to \infty}\frac{1}{N}\sum_{n=1}^N \mu(n)^2 f(x_n, x_{n+1}).
 $$
 where $\mu(\cdot)$ is the M\"obius function. 
 The key point is that $\eta:=(\mu(n)^2)$ is the characteristic function of the square free integers, which is $\nu$-generic for the associated Mirsky measure $\nu$ and the Mirsky measure is well described by its values on cylinders. 
 Let us state the following result, which is a particular case of Theorem \ref{thm:MA_Bfree}.
 In the following theorem, for simplicity, we will write $\mathbf{1}'B\mathbf{1}$
 by $\|B\|$, which is a norm for positive matrix $B$.

\begin{thm}
 The Lyapunov exponent $\psi(\beta)$ of the sequence of matrices $A_{\mu(n)^2}(\beta):=(e^{\beta \mu(n)^2 f(i, j)} )$ 
 is computed as follows
  $$
    \psi(\beta) = \left(1- 2F_1 -3F_{11} -4F_{111}\right) \log 2 + F_{1}\log \|A_1(\beta)\| + F_{11}\log \|A_1(\beta)^2\|+F_{111}\log \|A_1(\beta)^3\|,
 $$
 where $F_1,F_{11}, F_{111}$ are three numerical values \begin{eqnarray*}
     F_1&=& \frac{1}{\zeta_1(2)} -\frac{2}{\zeta_2(2)} + \frac{1}{\zeta_3(2)}=0.0881459,\\
     F_{11}&=& \frac{1}{\zeta_2(2)} - \frac{2}{\zeta_3(2)}=0.0716601, \\
    F_{111}&= &  \frac{1}{\zeta_3(2)}=0.125487,
 \end{eqnarray*}
 where  the functions $\zeta_a(s) $,  similar to the Riemann zeta function, are defined by 
  $$
 \zeta_a(s) =   \sum_{n=1}^\infty \frac{a^{\Omega(n)}}{n^{s}}  = \prod \Big(1- \frac{a}{p^{s}}\Big)^{-1}
 $$
 with $\Omega(n)$  being the number of primes with multiplicity contained in $n$.
The product here is taken over all prime numbers.
\end{thm}
  
 For the special case where $S=\{0,1\}$ and $f(x_0,x_1) = x_0x_1$, see Figure \ref{fig:0} for the graphs of  the Lyapunov exponent $\psi(\beta)$ as function of $\beta$ and of the multifractal spectrum $ \alpha\mapsto \dim \{x: \Phi(x)=\alpha\}$.  
   \begin{figure}[htb]
 	\centering
 	\includegraphics[width=0.45\linewidth]{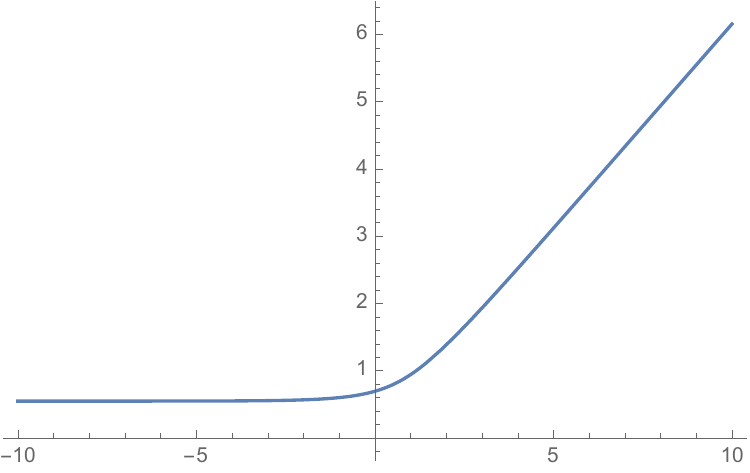}
 	\includegraphics[width=0.45\linewidth]{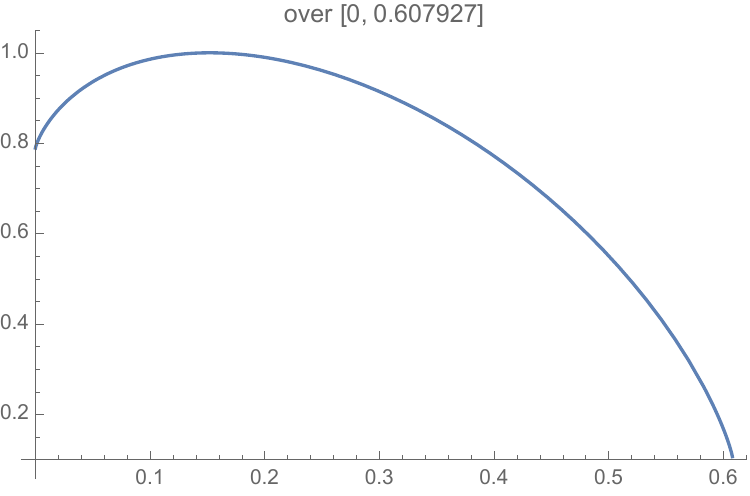}
 	\caption{ \small  Graphs of $\psi(\beta)$ and $\dim \{x: \Phi(x)=\alpha\}$ in the case: $f(x_0, x_1) = x_0x_1$ and $w_n=\mu(n)^2$.}
 	\label{fig:0}
 \end{figure}
 
 In this case, we have
$$
      A_0(\beta) = \begin{pmatrix}
      1 &1\\ 1&1
      \end{pmatrix} =\mathbf{1} \mathbf{1}', \qquad  
      A_1(\beta) = \begin{pmatrix}
      1 &1\\ 1&e^{\beta}
      \end{pmatrix}.
$$ 
Then
$$
  A_{11}(\beta) = \begin{pmatrix}
      2 &1+e^\beta\\ 1+e^\beta &1+ e^{2\beta}
      \end{pmatrix}, \qquad    A_{111}(\beta) = \begin{pmatrix}
      3+e^\beta &2+e^\beta + e^{2\beta}\\ 2+e^\beta + e^{2\beta} &1+ 2e^\beta +e^{3\beta}
      \end{pmatrix}.
$$
Thus
$$
  \|A_1(\beta)\|= 3+e^{\beta}, \quad
   \|A_{11}(\beta)\|= 5+2e^{\beta}+e^{2 \beta}, \quad  \|A_{111}(\beta)\|= 8+5e^{\beta} + 2 e^{2\beta} +e^{3\beta}
$$
On the other hand, since $\mathcal{R}=\{1, 11,111\}$, we have
$$
   \sum_{\mathbf{w}\in \mathcal{R}} (|\mathbf{w}|+1) F_{\mathbf{w}} = 2F_1 +3F_{11} +4F_{111}.
$$
By Theorem \ref{thm:MA_Bfree}, we get 
 $$
    \psi(\beta) = \left(1- 2F_1 -3F_{11} -4F_{111}\right) \log 2 + F_{1}\log \|A(\beta)\| + F_{11}\log \|A(\beta)^2\|+F_{111}\log \|A(\beta)^3\|.
 $$
 Consequently,
 $$
      \psi'(\beta) = F_{1} \frac{e^{\beta}}{3+e^{\beta} } + F_{11}\frac{2e^{\beta}+2 e^{2 \beta}}{5+2e^{\beta}+e^{2 \beta}}
      +F_{111}\frac{5e^{\beta} + 4 e^{2\beta} + 3e^{3\beta}}{8+5e^{\beta} + 2 e^{2\beta} +e^{3\beta}}.
 $$
 It is then clear that
 $$
    \psi'(-\infty)=0, \qquad \psi'(+\infty) = F_1+2F_{11}+3F_{111}=0.607927\cdots.
 $$
 See Figure \ref{fig:3} for the graph of $\psi'(\beta)$.

  \begin{figure}[htb]
 	\centering
 	\includegraphics[width=0.45\linewidth]{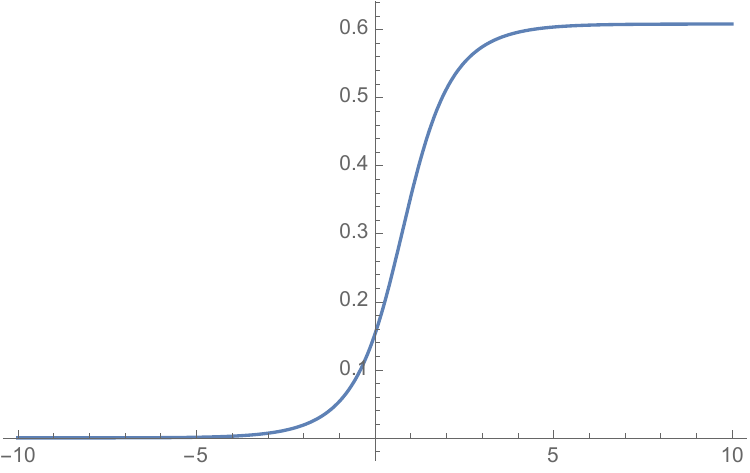}
 	\caption{ Graph of   $\psi'(\beta)$ in the case: $f(x_0, x_1) = x_0x_1$ and $w_n=\mu(n)^2$.}
 	\label{fig:3}
 \end{figure}

  The above multifractal analysis (equivalently the  computation of Lyapunov exponent $\psi$) for the weighted Birkhoff average  with weights $(\mu(n)^2)$
  is based on the fact that 
  the square of M\"{o}bius function $(\mu(n)^2)$ is generic.  Actually we only need the existence of  frequencies of 
  the patterns $1,11$ and $111$. 
  
  However  it is not known if the function of the M\"{o}bius function itself is generic for some invariant measure. 
 It is thus an open problem to determine the multifractal spectrum of
 the weighted Birkhoff average 
 $$
         \lim_{N\to \infty}\frac{1}{N}\sum_{n=1}^N \mu(n) f(x_n, x_{n+1}).
 $$ 
  In order to perform 
  the multifractal analysis for the   above Birkhoff average weighted by the M\"{o}bius function, we only need to know 
  the exact frequencies of the patterns $a, ab, abc$ with $a,b,c \in \{-1, 1\}$ (there are $14$ such patterns).  The frequencies of $-1$ and $1$ are known to be
  $\frac{3}{\pi^2}$. But it is an open problem to compute the frequencies of the patterns $ab, abc$ with $a,b,c \in \{-1, 1\}$ (cf. \cite{MRT2016}).

\end{document}